\documentclass[12pt]{article}
\usepackage{graphics}
\usepackage{graphicx}
\usepackage{color}
\usepackage[hypertex]{hyperref}
\usepackage{amsfonts}
\usepackage{eufrak}

\usepackage{amsmath}
\usepackage{amstext}
\usepackage{amsopn}
\usepackage{amsbsy}
\usepackage{amscd}
\usepackage{amsxtra}
\usepackage{amsthm}

\numberwithin{equation}{section}
\newtheorem{theorem}{Theorem}[section]
\newtheorem{lemma}[theorem]{Lemma}
\newtheorem{proposition}[theorem]{Proposition}
\newtheorem{corollary}[theorem]{Corollary}
\newtheorem{remark}[theorem]{Remark}

\newcommand{\BE}{\begin{equation}}
\newcommand{\BEN}{\begin{equation*}}
\newcommand{\EE}{\end{equation}}
\newcommand{\EEN}{\end{equation*}}
\newcommand{\BL}{\begin{lemma}}
\newcommand{\EL}{\end{lemma}}
\newcommand{\BT}{\begin{theorem}}
\newcommand{\ET}{\end{theorem}}
\newcommand{\BP}{\begin{proposition}}
\newcommand{\EP}{\end{proposition}}
\newcommand{\BC}{\begin{corollary}}
\newcommand{\EC}{\end{corollary}}
\newcommand{\BR}{\begin{remark}}
\newcommand{\ER}{\end{remark}}

\allowdisplaybreaks

\title{The virial theorem and ground state energy estimate of nonlinear Schr\"{o}dinger equations in $\mathbb{R}^2$ with square root and saturable nonlinearities in nonlinear optics }

\author{ Tai-Chia Lin\thanks{ Institute of Applied Mathematical Sciences and Center for Advanced Study in Theoretical Sciences (CASTS), National Taiwan University, Taipei, 10617, Taiwan. Email: tclin@math.ntu.edu.tw}
\ , Milivoj R. Beli\'c
\thanks{Texas A\&M University at Qatar, P.O. Box 23874, Doha,
Qatar. }\ , Milan S. Petrovi\'c
\thanks{Institute of Physics, P.O. Box 57, 11001 Belgrade,
Serbia}\ ,
Hichem Hajaiej
\thanks{New York University－Shanghai, Pudong New District, Shanghai 200122, China. } \ , \\
Goong Chen$^\dagger$\thanks{Department of Mathematics and Institute for Quantum Science and Engineering, Texas A\&M University,
College Station, Texas 77843, USA. } }

\date{\today}

\begin{document}

\maketitle

\begin{abstract}
\noindent The virial theorem is a nice property for the linear Schr\"{o}dinger equation in atomic and molecular physics as it gives an elegant ratio between the kinetic and potential energies and is useful in assessing the quality of numerically computed eigenvalues. If the governing equation is a nonlinear Schr\"{o}dinger equation with power-law nonlinearity, then a similar ratio can be obtained but there seems no way of getting any eigenvalue estimate. It is surprising as far as we are concerned that when the nonlinearity is either square-root or saturable nonlinearity (not a power-law), one can develop a virial theorem and eigenvalue estimate of nonlinear Schr\"{o}dinger (NLS) equations in ${{\mathbb{R}}^{2}}$ with square-root and saturable nonlinearity, respectively. Furthermore, we show here that the eigenvalue estimate can be used to obtain the 2nd order term (which is of order $\ln\Gamma$) of the lower bound of the ground state energy as the coefficient $\Gamma$ of the nonlinear term tends to infinity.
\end{abstract}

\section{Introduction}\label{s1}
\medskip
\noindent

The nonlinear Schr\"{o}dinger (NLS) equation, a nonlinear variation of the Schr\"{o}dinger equation, is a universal model in nonlinear science and mathematics. Such an equation can be represented as follows:
\BE\label{ns-1}
i\frac{\partial A}{\partial z}+ \Delta A+\Gamma f\left( {{\left| A \right|}^{2}} \right)A=0\,,
\EE
for $z>0,\text{ }x=\left( {{x}_{1}},\cdots ,{{x}_{d}} \right)\in {{\mathbb{R}}^{d}},\text{ }$ where $A=A\left( x,z \right)\in \mathbb{C}$,$\Gamma \in \mathbb{R}$, $d\geq 1$, $\Delta =\sum\limits_{j=1}^{d}{\partial _{{{x}_{j}}}^{2}}$ and the function $f$ denotes the nonlinearity. Physically, $A$ is the wave function, $d$ is the transverse dimension and $\Gamma $ is the strength of nonlinearity. Here we study the case of $d=2$ for equation (\ref{ns-1}) which is non-integrable. But note that in the case of $d=1$, equation (\ref{ns-1}) becomes integrable, and thus can be investigated by different methods of the inverse scattering theory~\cite{ablowitz}. For the nonlinearity in function $f$, we consider the square-root and saturable nonlinearities in the following forms:
\begin{enumerate}
\item[(I)]~~{\bf square-root nonlinearity:}~$f\left( s \right)=1-\frac{1}{\sqrt{1+s}}$ for $s>0$, \\
\item[(II)]~{\bf saturable nonlinearity:}~$f\left( s \right)=1-\frac{1}{1+s}$ for $s>0$,
\end{enumerate}
\medskip
which describe {\em narrow-gap semiconductors} (\cite{PPB-oc07, S-prb2000}) and {\em photorefractive media} (\cite{GH-josab97, K-pra92, K-prl65, MMP-cvpde12, MD-prl68, MMSN-josab07, PBDK-lpor_2011}), respectively. Equation (\ref{ns-1}) can be represented as $i\frac{\partial A}{\partial z}=\frac{\delta E\left[ A \right]}{\delta A}$, where
$E\left[ A \right]=\frac{1}{2}\int_{{{\mathbb{R}}^{2}}}{{{\left| \nabla A \right|}^{2}}-\Gamma F\left( {{\left| A \right|}^{2}} \right)dx}$
and $F\left( I \right)=\int_{0}^{I}{f\left( s \right)ds}$. Besides, the total energy $E=K+P$ can be denoted as the sum of the kinetic energy $K$ and the potential energy $P$, where the kinetic energy is
\BE\label{ke1}
K\left[ A \right]=\frac{1}{2}\int_{{{\mathbb{R}}^{2}}}{{{\left| \nabla A \right|}^{2}}dx}\,,
\EE
and the potential energy is
\BE\label{pe1}
P\left[ A \right]=-\frac{\Gamma }{2}\int_{{{\mathbb{R}}^{2}}}{F\left( {{\left| A \right|}^{2}} \right)dx}\,.
\EE

To see solitons of equation (\ref {ns-1}), we may set $A\left( x,z \right)={{e}^{i\lambda z}}u\left( x \right)$   for $x\in {{\mathbb{R}}^{2}}$ and $z>0$, where  $\lambda \in \mathbb{R}$ is a constant and  $u=u\left( x \right)$ is a real-valued function. Then by (\ref{ns-1}), we get the following nonlinear eigenvalue problem:
\BE\label{evp1}
\Delta u+\Gamma f\left( {{u}^{2}} \right)u=\lambda u \quad\hbox{ in }\: {{\mathbb{R}}^{2}}\,,
\EE
where $\lambda $ is an eigenvalue and $u$ is the associated eigenfunction. When function $f\left( s \right)={{s}^{\frac{p-1}{2}}}$, $p>1$ is of power law nonlinearity~\cite{CL-cmp82, W-cmp83, W-siamjma85}, the eigenvalue $\lambda $ can be apriori chosen as a positive number because we may set $u\left( x \right)={{\left( \frac{\lambda }{\Gamma } \right)}^{\frac{1}{p-1}}}U\left( \sqrt{\lambda }x \right)$ and transform equation (\ref{evp1}) into $\Delta U-U+{{U}^{p}}=0$ in ${{\mathbb{R}}^{2}}$, which has a unique positive solution $U$. However, when function $f$ is of square-root and saturable nonlinearity, the eigenvalue $\lambda $ cannot be any positive number. One na\"{i}ve counterexample is to set $\lambda =\Gamma>0 $ and equation (\ref{evp1}) has only zero solution because of the standard Liouville theorem. This motivates us to study the estimate of the eigenvalue $\lambda$ of the ground state of equation (\ref{evp1}) with the square-root and saturable nonlinearity of function~$f$.

The {\em virial theorem} of linear Schr\"{o}dinger equations can be formulated as the ratio of the kinetic energy and the potential energy of linear Schrodinger equations. Such a theorem plays a useful role in assessing the quality of numerical solutions of the eigenvalues of linear Schrodinger equations which is important in quantum mechanics (cf.~\cite{CDH-v-2014}). To develop a virial theorem for the nonlinear eigenvalue problem (\ref{evp1}), we consider $\frac{K}{P}$ the ratio of the kinetic energy and the potential energy defined in (\ref{ke1}) and (\ref{pe1}), respectively.  For the power-law nonlinearity of function $f\left( s \right)={{s}^{\frac{p-1}{2}}}$, $p>1$, the ratio is denoted as
\[a=\frac{\frac{1}{2}\int_{{{\mathbb{R}}^{2}}}{{{\left| \nabla u \right|}^{2}}}}{-\frac{\Gamma }{p+1}\int_{{{\mathbb{R}}^{2}}}{{{\left| u \right|}^{p+1}}}}\,,\]
where $u\in {{H}^{1}}\left( {{\mathbb{R}}^{2}} \right)\cap {{L}^{p+1}}\left( {{\mathbb{R}}^{2}} \right)$ satisfies
\BE\label{pwl-nls}
\Delta u+\Gamma {{\left| u \right|}^{p-1}}u=\lambda u\quad\hbox{ in }\quad {{\mathbb{R}}^{2}}\,.
\EE
Then by direct calculation, we get
\BE\label{pwl-pid1}
\frac{2}{p+1}\Gamma \int_{{{\mathbb{R}}^{2}}}{{{\left| u \right|}^{p+1}}}=\lambda \int_{{{\mathbb{R}}^{2}}}{{{u}^{2}}}\,,
\EE
and
\BE\label{pwl-id2}
-\int_{{{\mathbb{R}}^{2}}}{{{\left| \nabla u \right|}^{2}}}+\Gamma \int_{{{\mathbb{R}}^{2}}}{{{\left| u \right|}^{p+1}}}=\lambda \int_{{{\mathbb{R}}^{2}}}{{{u}^{2}}}\,.
\EE
Here (\ref{pwl-pid1}) is the Pohozaev identity of (\ref{pwl-nls}) (cf.~\cite{P-smd1965}), and (\ref{pwl-id2}) is the constraint of the Nehari manifold of (\ref{pwl-nls}) (cf.~\cite{BL-arma83-1, BL-arma83-2}). Combining (\ref{pwl-pid1}) and (\ref{pwl-id2}), we have
\[\int_{{{\mathbb{R}}^{2}}}{{{\left| \nabla u \right|}^{2}}}=\Gamma \left( 1-\frac{2}{p+1} \right)\int_{{{\mathbb{R}}^{2}}}{{{\left| u \right|}^{p+1}}}\,,\]
implying that the ratio is $a=\frac{1}{2}\left( 1-p \right)$ so we may represent the virial relation of (\ref{pwl-nls}) as $\frac{K}{P}=a=\frac{1}{2}\left( 1-p \right)$. When $1<p<3$ (which is of the subcritical case and the existence of ground state is proved in~\cite{CL-cmp82, W-cmp83}), the ratio $a$ is located on the interval $(-1,0)$, which is the same interval (up to boundary points) as the virial theorem of linear Schr\"{o}dinger equations with Coulomb potentials (see Section~4 of \cite{CDH-v-2014}). However, it seems impossible to get any estimate of the eigenvalue $\lambda$ from identity (\ref{pwl-pid1}) with $p>1$. This stimulates us to study the different types of nonlinearities such as the square-root and saturable types here.

For the square-root~$(I)$ and saturable~$(II)$ nonlinearities of function $f$, we define the ratio as follows
\BE\label{af-1}
\alpha =\frac{\int_{{{\mathbb{R}}^{2}}}{{{\left| \nabla u \right|}^{2}}}}{-\Gamma \int_{{{\mathbb{R}}^{2}}}{{{\left( \sqrt{1+{{u}^{2}}}-1 \right)}^{2}}}}\,,
\EE
and
\BE\label{ba-1}
\beta =\frac{\int_{{{\mathbb{R}}^{2}}}{{{\left| \nabla u \right|}^{2}}}}{-\Gamma \int_{{{\mathbb{R}}^{2}}}{\left[ {{u}^{2}}-\ln \left( 1+{{u}^{2}} \right) \right]}}\,,
\EE
respectively. Note that the ratio $\alpha$ is for the eigenvalue problem (with square-root nonlinearity)
\BE\label{eq-a.1}
\Delta u+\Gamma \left( 1-\frac{1}{\sqrt{1+{{u}^{2}}}} \right)u=\tilde{\lambda }\,u \quad\hbox{ in }\quad {{\mathbb{R}}^{2}}\,,
\EE
and the ratio $\beta $ is for the eigenvalue problem (with saturable nonlinearity)
\BE\label{eq-1.1}
\Delta u+\Gamma \left( 1-\frac{1}{1+{{u}^{2}}} \right)u=\hat{\lambda }\,u \quad\hbox{ in }\quad {{\mathbb{R}}^{2}}\,,
\EE
where $\tilde{\lambda }$ and $\hat{\lambda }$ are the respective eigenvalues, and $u\in {{H}^{1}}\left( {{\mathbb{R}}^{2}} \right)$ is the associated eigenfunction. In this paper, we first prove that the ratios $\alpha $ defined in (\ref{af-1}) and $\beta $ defined in (\ref{ba-1}) must be located in the interval $(-1,0)$, which is the same interval (up to boundary points) as the virial theorem of linear Schr\"{o}dinger equations with suitable potentials (see Section~4 of \cite{CDH-v-2014}). Then we use the ratios $\alpha$ and $\beta$ to derive the eigenvalue estimate of the ground states of (\ref{eq-a.1}) and (\ref{eq-1.1}), respectively (see Theorems~\ref{theomII} and~\ref{theomI}).

For the estimates of the ratio $\alpha$ and the ground state eigenvalue $\tilde{\lambda }$, the results are stated as follows.
\BT\label{theomII}
Let $u\in {{H}^{1}}\left( {{\mathbb{R}}^{2}} \right)$ be an eigenfunction and $\tilde{\lambda }\in \mathbb{R}$ be the eigenvalue of problem (\ref{eq-a.1}). Suppose $\Gamma >0$.
Then the ratio $\alpha =\frac{\int_{{{\mathbb{R}}^{2}}}{{{\left| \nabla u \right|}^{2}}}}{-\Gamma \int_{{{\mathbb{R}}^{2}}}{{{\left( \sqrt{1+{{u}^{2}}}-1 \right)}^{2}}}}$ defined in (\ref{af-1}) satisfies
\BE\label{a.3}
-1<\alpha <0\,.
\EE
Moreover, the eigenvalue $\tilde{\lambda }$ has the following estimate:
\begin{enumerate}
\item[(i)]~~~If $-\frac{1}{2}\le \alpha <0$, then $0<\tilde{\lambda }\le \Gamma \left( 1+\alpha  \right)$,
\item[(ii)]~~If $-1<\alpha <-\frac{1}{2}$, then $0<\tilde{\lambda }\le \Gamma \left[ 1+\alpha +{{\left( -\frac{1+2\alpha }{3} \right)}^{\frac{3}{2}}} \right]$,
\item[(iii)]~If $-1<\alpha <-\frac{1}{2}$ and $\left\| u \right\|_{\infty }^{2}\le {{\left( 1+2\alpha  \right)}^{-2}}-1$, then $0<\tilde{\lambda }\le \Gamma \left( 1+\alpha  \right)$, where ${{\left\| u \right\|}_{\infty }}=\underset{x\in {{\mathbb{R}}^{2}}}{\mathop{\max }}\,\left| u\left( x \right) \right|$.
\end{enumerate}
\ET
\noindent For the estimates of the ratio $\beta $ and the ground state eigenvalue $\hat{\lambda }$, the results are stated as follows.
\BT\label{theomI}
Let $u\in {{H}^{1}}\left( {{\mathbb{R}}^{2}} \right)$ be the eigenfunction and $\hat{\lambda }\in \mathbb{R}$ be the eigenvalue of problem (\ref{eq-1.1}). Suppose $\Gamma >0$.
Then the ratio $\beta =\frac{\int_{{{\mathbb{R}}^{2}}}{{{\left| \nabla u \right|}^{2}}}}{-\Gamma \int_{{{\mathbb{R}}^{2}}}{\left[ {{u}^{2}}-\ln \left( 1+{{u}^{2}} \right) \right]}}$ and the eigenvalue $\hat{\lambda }$ satisfy
\BE\label{1.3}
-1<\beta <0\,,
\EE
\BE\label{1.4}
0<\hat{\lambda }\le {{\left( 1+\frac{1}{2}\beta  \right)}^{2}}\Gamma\,.
\EE
\ET
\BR
From Theorems~\ref{theomII} and \ref{theomI}, the ratios $\alpha$ and $\beta$ satisfy $-1< \alpha, \beta < 0$ and is located on the same interval (up to boundary points) as the virial theorem of linear Schr\"{o}dinger equations with Coulomb potentials (see Section~4 of \cite{CDH-v-2014}).
\ER

With saturable nonlinearity, the ground state of the nonlinear Schr\"{o}dinger equation (\ref {ns-1}) is defined as the minimizer of the following problem:
\[\text{Minimize}\left\{ E\left[ u \right]:u\in {{H}^{1}}\left( {{\mathbb{R}}^{2}} \right), \|u\|_2=1 \right\}\,, \]
and the ground state existence can be proved using an energy estimate method~\cite{LC-jmp-2014},
where $E\left[ u \right]=\frac{1}{2}\int_{{{\mathbb{R}}^{2}}}{{{\left| \nabla u \right|}^{2}}-\Gamma\,\left[ {{u}^{2}}-\ln \left( 1+{{u}^{2}} \right) \right] dx}$ and $\| u \|_2^2=\int_{{{\mathbb{R}}^{2}}}{{{u}^{2}}dx}$. Such a ground state satisfies equation (\ref{eq-1.1}) and the eigenvalue $\hat{\lambda }$ comes from the Lagrange multiplier corresponding to the ${{L}^{2}}$-norm constraint ${{\left\| u \right\|}_{2}}=1$. The result of ground state existence is stated as follows.
\newline
\noindent {\bf Theorem~A.}(cf.~\cite{LC-jmp-2014})~~Consider the following minimization problem:
\BE\label{gs-sat}
{{e }_{\Gamma }}=\text{Minimize}\left\{ {{E}}\left[ u \right]:u\in {{H}^{1}}\left( {{\mathbb{R}}^{2}} \right),{{\left\| u \right\|}_{2}}=1 \right\}\,,
\EE
where ${{E}}\left[ u \right]=\int_{{{\mathbb{R}}^{2}}}{{{\left| \nabla u \right|}^{2}}-\Gamma \left[ {{u}^{2}}-\ln \left( 1+{{u}^{2}} \right) \right].}$  Let ${{T}_{1}}$ be the following positive constant:
\[{{T}_{1}}=\underset{\begin{matrix}
   w\in {{H}^{1}}\left( {{\mathbb{R}}^{2}} \right)  \\
   {{\left\| w \right\|}_{2}}=1  \\
\end{matrix}}{\inf}\,\frac{\int_{{{\mathbb{R}}^{2}}}{{{\left| \nabla w \right|}^{2}}}}{\int_{{{\mathbb{R}}^{2}}}{\left[ {{w}^{2}}-\ln\left( 1+{{w}^{2}} \right) \right]}}\,.\]
Then
\begin{enumerate}
\item[(i)]~~If $\Gamma < {{T}_{1}}$, then ${{e }_{\Gamma }}=0$ can not be attained by a minimizer, i.e., problem (\ref{gs-sat}) has no ground state.
\item[(ii)]~If $\Gamma > {{T}_{1}}$, then ${{e }_{\Gamma }}<0$ and there exists a minimizer of (\ref{gs-sat}) denoted as $U=U(r)$ which is radially symmetric and monotone decreasing for $r>0$.
\end{enumerate}
\medskip

\noindent The positivity of ${{T}_{1}}$ comes from the fact that ${{w}^{2}}-\ln\left( 1+{{w}^{2}} \right)\le \tfrac{1}{2}{{w}^{4}}$ for $w\in \mathbb{R}$ and the Gagliardo-Nirenberg inequality (cf.~\cite{G-rm58, N-cpam55}):
\BE\label{NG-est}
\int_{{{\mathbb{R}}^{2}}}{{{\left| \nabla w \right|}^{2}}}\ge {{C}_{0}}\left\| w \right\|_{2}^{-2}\left\| w \right\|_{4}^{4}\quad\text{ for }\quad w\in {{H}^{1}}\left( {{\mathbb{R}}^{2}} \right)\,,
\EE
for a positive constant ${{C}_{0}}$. Hereafter, the norm ${{\left\| \cdot  \right\|}_{p}}$ denotes ${{\left\| w \right\|}_{p}}={\left( \int_{{{\mathbb{R}}^{2}}}{{{\left| w \right|}^{p}}} \right)}^{1/p}$ for $p>1$. Using the eigenvalue estimate of Theorem~\ref{theomI}, we can derive the ground state energy estimate of ${{e}_{\Gamma }}$ as follows.
\BT\label{theomIII}
Let ${{e}_{\Gamma }}$ be the ground state energy defined in (\ref{gs-sat}). Then as $\Gamma \to \infty $, ${{e}_{\Gamma }}=-\frac{\Gamma }{2}\left( 1+{{o}_{\Gamma }}\left( 1 \right) \right)$, where ${{o}_{\Gamma }}\left( 1 \right)$ is a small quantity tending to zero as $\Gamma $ goes to infinity. Furthermore, for $\sigma \in \left( 0,1 \right)$, there exists a positive constant ${{\Gamma }_{\sigma }}$ such that
\BE\label{lb-est1}
{{e}_{\Gamma }}\ge -\frac{\Gamma }{2}+\sigma \frac{{{T}_{1}}}{2}\ln \Gamma +{{C}_{0}}\quad\hbox{ for }\quad \Gamma >{{\Gamma }_{\sigma }}\,,
\EE
where ${{C}_{0}}$ is a constant independent of $\Gamma $.
\ET
\medskip
\noindent Note that in (\ref{lb-est1}), the 2nd order term of the lower bound of ${{e}_{\Gamma }}$ is of order $\ln \Gamma $ which goes to positive infinity as $\Gamma $ tends to infinity.

With a square-root nonlinearity, the ground state of the nonlinear Schr\"{o}dinger equation (\ref {ns-1}) is defined as the minimizer of the following problem:
\[\text{Minimize}\left\{ E\left[ u \right]:u\in {{H}^{1}}\left( {{\mathbb{R}}^{2}} \right), \|u\|_2=1 \right\}\,,\]
where $E\left[ u \right]=\frac{1}{2}\int_{{{\mathbb{R}}^{2}}}{{{\left| \nabla u \right|}^{2}}-\Gamma {{\left( \sqrt{1+{{u}^{2}}}-1 \right)}^{2}}dx}$. Such a ground state satisfies equation (\ref{eq-a.1}) and the eigenvalue $\tilde{\lambda }$ comes from the Lagrange multiplier corresponding to the ${{L}^{2}}$-norm constraint ${{\left\| u \right\|}_{2}}=1$. We may generalize the argument of~\cite{LC-jmp-2014} to prove the existence of the ground state and obtain the following result.
\newline
\noindent {\bf Theorem~B.}~~ Consider the following minimization problem:
\BE\label{gs-2}
{{\tilde{e}}_{\Gamma }}=\text{Minimize}\left\{ {{{E}}}\left[ u \right]:u\in {{H}^{1}}\left( {{\mathbb{R}}^{2}} \right)\,, \| u \|_2=1 \right\}\,,
\EE
where $E\left[ u \right]=\frac{1}{2}\int_{{{\mathbb{R}}^{2}}}{{{\left| \nabla u \right|}^{2}}-\Gamma {{\left( \sqrt{1+{{u}^{2}}}-1 \right)}^{2}}dx}$.
Let ${{T}_{2}}$ be the following positive constant:
$$
{{T}_{2}}=\underset{\begin{matrix}
   w\in {{H}^{1}}\left( {{\mathbb{R}}^{2}} \right)  \\
   {{\left\| w \right\|}_{2}}=1  \\
\end{matrix}}{\mathop{\inf }}\,\frac{\int_{{{\mathbb{R}}^{2}}}{{{\left| \nabla w \right|}^{2}}}}{\int_{{{\mathbb{R}}^{2}}}{{{\left( \sqrt{1+{{w}^{2}}}-1 \right)}^{2}}}}\,.
$$
\begin{enumerate}
\item[(i)]~~If $\Gamma < {{T}_{2}}$, then ${{\tilde{e}}_{\Gamma }}=0$ can not be attained by a minimizer, i.e., problem (\ref{gs-2}) has no ground state.
\item[(ii)]~If $\Gamma > {{T}_{2}}$, then ${{\tilde{e}}_{\Gamma }}<0$ and there exists a minimizer $U=U(r)$ of (\ref{gs-2}) which is radially symmetric and monotone decreasing for $r>0$.
\end{enumerate}

\noindent Here the positivity of ${{T}_{2}}$ comes from (\ref{NG-est}), i.e., the Gagliardo-Nirenberg inequality (cf.~\cite{G-rm58, N-cpam55}) and the fact that ${{\left( \sqrt{1+{{w}^{2}}}-1 \right)}^{2}}\le \frac{1}{4}w^4$ for $w\in \mathbb{R}$. The proof of Theorem~B is similar to that of Theorem~A so we need only provide a brief sketch of the proof in Appendix~II. We can use the eigenvalue estimate of Theorem~\ref{theomII} to derive the ground state energy estimate of ${{\tilde{e}}_{\Gamma }}$ as follows.
\BT\label{theomIV}
Let ${{\tilde{e}}_{\Gamma }}$ be the ground state energy defined in (\ref{gs-2}). Then as $\Gamma \to \infty $, ${{\tilde{e}}_{\Gamma }}=-\frac{\Gamma }{2}\left( 1+{{o}_{\Gamma }}\left( 1 \right) \right)$, where ${{o}_{\Gamma }}\left( 1 \right)$ is a small quantity tending to zero as $\Gamma $ goes to infinity. Furthermore, there exists a positive constant ${{\Gamma }_{0}}$ such that
\BE\label{lb-est2}
{{\tilde{e}}_{\Gamma }}\ge -\frac{\Gamma }{2}+\frac{{{T}_{2}}}{2}\ln \Gamma +{{C}_{0}}\quad\hbox{ for }\quad \Gamma >{{\Gamma }_{0}}\,,
\EE
where ${{C}_{0}}$ is a constant independent of $\Gamma $.
\ET
\medskip
\noindent Note that in (\ref{lb-est2}), the 2nd order term of the lower bound of ${{\tilde{e}}_{\Gamma }}$ is also of order $\ln \Gamma$ (same as that of ${{e}_{\Gamma }}$) which goes to positive infinity as $\Gamma $ tends to infinity. On the other hand, the difference between the ground state energy estimate (\ref{lb-est1}) and (\ref{lb-est2}) comes from that of the eigenvalue estimate (\ref{1.4}) (see Theorem~\ref{theomI}) and Theorem~\ref{theomII}~(i).

The rest of this paper is organized as follows: The proofs of Theorem~\ref{theomII} and~\ref{theomIV} are given in Sections~\ref{sqrt1} and ~\ref{sqrt2}, respectively. We provide the proofs of Theorem~\ref{theomI} and~\ref{theomIII} in Sections~\ref{saue1} and~\ref{saue2}, respectively. Brief concluding remarks are given in Section~\ref{CR}.

\section{Proof of Theorem~\ref{theomII}}\label{sqrt1}
We multiply (\ref{eq-a.1}) by $u$ and integrate it over ${{\mathbb{R}}^{2}}$. Then using integration by parts, we get
\BE\label{a.6}
-\int_{{{\mathbb{R}}^{2}}}{{{\left| \nabla u \right|}^{2}}}+\Gamma \int_{{{\mathbb{R}}^{2}}}{\left( 1-\frac{1}{\sqrt{1+{{u}^{2}}}} \right)}{{u}^{2}}=\tilde{\lambda }\left\| u \right\|_{2}^{2}\,,
\EE
where $\left\| u \right\|_{2}^{2}=\int_{{{\mathbb{R}}^{2}}}{{{u}^{2}}}$.
On the other hand, we may multiply (\ref{eq-a.1}) by $x\cdot \nabla u$ and integrate it over ${{\mathbb{R}}^{2}}$, where $x\cdot \nabla u=\sum\limits_{j=1}^{2}{{{x}_{j}}{{\partial }_{j}}u}$ and ${{\partial }_{j}}u=\frac{\partial u}{\partial {{x}_{j}}}$. Then using integration by parts, we can derive the Pohozaev identity as follows
\BE\label{a.7}
\tilde{\lambda }\left\| u \right\|_{2}^{2}=\Gamma \int_{{{\mathbb{R}}^{2}}}{{{\left( \sqrt{1+{{u}^{2}}}-1 \right)}^{2}}}\,.
\EE
Combining (\ref{a.6}) and (\ref{a.7}), we have
\[-\int_{{{\mathbb{R}}^{2}}}{{{\left| \nabla u \right|}^{2}}}+\Gamma \int_{{{\mathbb{R}}^{2}}}{\left( 1-\frac{1}{\sqrt{1+{{u}^{2}}}} \right)}{{u}^{2}}=\Gamma {{\int_{{{\mathbb{R}}^{2}}}{\left( \sqrt{1+{{u}^{2}}}-1 \right)}}^{2}}\,,\]
which implies
\BE\label{a.8}
\int_{{{\mathbb{R}}^{2}}}{{{\left| \nabla u \right|}^{2}}}=\Gamma \int_{{{\mathbb{R}}^{2}}}{\left( 1-\frac{1}{\sqrt{1+{{u}^{2}}}} \right)}{{u}^{2}}-\Gamma \int_{{{\mathbb{R}}^{2}}}{{{\left( \sqrt{1+{{u}^{2}}}-1 \right)}^{2}}}\,.
\EE
Hence
\BE\label{a.9}
\begin{array}{ll}
\alpha &=\frac{\int_{{{\mathbb{R}}^{2}}}{{{\left| \nabla u \right|}^{2}}}}{-\Gamma \int_{{{\mathbb{R}}^{2}}}{{{\left( \sqrt{1+{{u}^{2}}}-1 \right)}^{2}}}} \\
 & =1-\frac{\int_{{{\mathbb{R}}^{2}}}{\left( 1-\frac{1}{\sqrt{1+{{u}^{2}}}} \right){{u}^{2}}}}{\int_{{{\mathbb{R}}^{2}}}{{{\left( \sqrt{1+{{u}^{2}}}-1 \right)}^{2}}}}
\end{array}
\EE
It is obvious that $\alpha <0$ because $u$ is nontrivial.

To prove $\alpha >-1$, we define a function $f\left( s \right)=\frac{\left( 1-\frac{1}{\sqrt{1+s}} \right)s}{{{\left( \sqrt{1+s}-1 \right)}^{2}}}$ for $s>0$. Then by direct calculation, we have $f\left( s \right)=\frac{s}{1+s-\sqrt{1+s}}$ and ${f}'\left( s \right)={{\left( 1+s-\sqrt{1+s} \right)}^{-2}}\frac{1}{\sqrt{1+s}}\left( \sqrt{1+s}-1-\frac{1}{2}s \right)<0$ for $s>0$, which gives ${f}'\left( s \right)<0$ for $s>0$. Here we have used the fact that $\sqrt{1+s}<1+\frac{1}{2}s$ for $s>0$. Besides, $\underset{s\to 0+}{\mathop{\lim }}\,f\left( s \right)=2$ and $\underset{s\to \infty }{\mathop{\lim }}\,f\left( s \right)=1$ are trivial by direct calculation. Consequently, $f\left( s \right)<2$ for $s>0$ and
\[\int_{{{\mathbb{R}}^{2}}}{\left( 1-\frac{1}{\sqrt{1+{{u}^{2}}}} \right)}{{u}^{2}}=\int_{{{\mathbb{R}}^{2}}}{{{\left( \sqrt{1+{{u}^{2}}}-1 \right)}^{2}}}f\left( {{u}^{2}} \right)\le 2\int_{{{\mathbb{R}}^{2}}}{{{\left( \sqrt{1+{{u}^{2}}}-1 \right)}^{2}}}\]
which gives $\alpha \ge -1$. Here we have used (\ref{a.9}).

Now we prove $\alpha >-1$ by contradiction. Suppose $\alpha =-1$, i.e.,
\BE\label{a.10}
\int_{{{\mathbb{R}}^{2}}}{\left( 1-\frac{1}{\sqrt{1+{{u}^{2}}}} \right)}{{u}^{2}}=2\int_{{{\mathbb{R}}^{2}}}{{{\left( \sqrt{1+{{u}^{2}}}-1 \right)}^{2}}}\,.
\EE
Due to $u\not\equiv 0$, the unique continuation theorem of equation (\ref{eq-a.1}) implies that there exists a ball ${{B}_{0}}={{B}_{{{R}_{0}}}}\left( {{x}_{0}} \right)$ with radius ${{R}_{0}}>0$ and center ${{x}_{0}}\in {{\mathbb{R}}^{2}}$ such that ${{u}^{2}}\left( x \right)>0$ for $x\in {{B}_{0}}$. Then there exist ${{\varepsilon }_{0}}>0$ and a smaller ball ${{B}_{1}}={{B}_{{{R}_{1}}}}\left( {{x}_{0}} \right)\subset \subset {{B}_{0}}$ such that ${{u}^{2}}\left( x \right)\ge {{\varepsilon }_{0}}>0$ for $x\in {{B}_{1}}$. Because of ${f}'\left( s \right)<0$ for $s>0$, $\underset{s\to 0+}{\mathop{\lim }}\,f\left( s \right)=2$ and $\underset{s\to \infty }{\mathop{\lim }}\,f\left( s \right)=1$, there exists ${{\delta }_{0}}>0$ such that $f\left( {{u}^{2}}\left( x \right) \right)\le 2-{{\delta }_{0}}$ for $x\in {{B}_{1}}$. Hence (\ref{a.10}) causes the following contradiction:
\begin{eqnarray*}
2\int_{{{\mathbb{R}}^{2}}}{{{\left( \sqrt{1+{{u}^{2}}}-1 \right)}^{2}}}&=& \int_{{{\mathbb{R}}^{2}}}{\left( 1-\frac{1}{\sqrt{1+{{u}^{2}}}} \right)}{{u}^{2}} \\
 &=& \int_{{{B}_{1}}}{\left( 1-\frac{1}{\sqrt{1+{{u}^{2}}}} \right)}{{u}^{2}}+\int_{B_{1}^{c}}{\left( 1-\frac{1}{\sqrt{1+{{u}^{2}}}} \right)}{{u}^{2}} \\
 &=& \int_{{{B}_{1}}}{f\left( {{u}^{2}} \right)}{{\left( \sqrt{1+{{u}^{2}}}-1 \right)}^{2}}+\int_{B_{1}^{c}}{f\left( {{u}^{2}} \right)}{{\left( \sqrt{1+{{u}^{2}}}-1 \right)}^{2}} \\
 &\le& \left( 2-{{\delta }_{0}} \right)\int_{{{B}_{1}}}{{{\left( \sqrt{1+{{u}^{2}}}-1 \right)}^{2}}}+2\int_{B_{1}^{c}}{{{\left( \sqrt{1+{{u}^{2}}}-1 \right)}^{2}}} \\
 &=& 2\int_{{{\mathbb{R}}^{2}}}{{{\left( \sqrt{1+{{u}^{2}}}-1 \right)}^{2}}}-{{\delta }_{0}}\int_{{{B}_{1}}}{{{\left( \sqrt{1+{{u}^{2}}}-1 \right)}^{2}}} \\
 &<& 2\int_{{{\mathbb{R}}^{2}}}{{{\left( \sqrt{1+{{u}^{2}}}-1 \right)}^{2}}}\,,
\end{eqnarray*}
where $B_{1}^{c}={{\mathbb{R}}^{2}}-{{B}_{1}}$ is the complement of ${{B}_{1}}$. Therefore, we have completed the proof for the case $-1<\alpha <0$, i.e., (\ref{a.3}).

To prove Theorem~\ref{theomII}~(i) and~(iii), we substitute (\ref{af-1}) into (\ref{a.6}) and get
\BE\label{a.11}
\tilde{\lambda }\left\| u \right\|_{2}^{2}=\alpha \Gamma \int_{{{\mathbb{R}}^{2}}}{{{\left( \sqrt{1+{{u}^{2}}}-1 \right)}^{2}}}+\Gamma \int_{{{\mathbb{R}}^{2}}}{\left( 1-\frac{1}{\sqrt{1+{{u}^{2}}}} \right)}{{u}^{2}}=\Gamma \int_{{{\mathbb{R}}^{2}}}{g\left( {{u}^{2}} \right)}\,,
\EE
where $g\left( s \right)=\alpha {{\left( \sqrt{1+s}-1 \right)}^{2}}+s-\frac{s}{\sqrt{1+s}}$  for $s>0$. Let $h\left( s \right)=g\left( s \right)-\left( 1+\alpha  \right)s$ for $s>0$. Then $h\left( s \right)=2\alpha \left( 1-\sqrt{1+s} \right)-\frac{s}{\sqrt{1+s}}$ for $s>0$. By a direct calculation, we get ${h}'\left( s \right)=\frac{-1}{2\sqrt{1+s}}\left( 1+2\alpha +\frac{1}{1+s} \right)$ for $s>0$. Suppose $-\frac{1}{2}\le \alpha <0$. Then ${h}'\left( s \right)<0$ for $s>0$, implying $h\left( s \right)<h\left( 0 \right)=0$, i.e. $g\left( s \right)<\left( 1+\alpha  \right)s$ for $s>0$, which can be used in (\ref{a.11}) to complete the proof of Theorem~\ref{theomII}~(i). On the other hand, suppose $-1<\alpha <-\frac{1}{2}$. Then $h\left( s \right)\le 0$ for $0\le s\le {{s}_{\alpha }}$, and $h\left( s \right)>0$ for $s>{{s}_{\alpha }}$, where ${{s}_{\alpha }}>0$ satisfies $2\alpha \left[ \sqrt{1+{{s}_{\alpha }}}-\left( 1+{{s}_{\alpha }} \right) \right]-{{s}_{\alpha }}=0$, i.e. ${{s}_{\alpha }}=\frac{1}{{{\left( 1+2\alpha  \right)}^{2}}}-1>0$. Consequently, $h\left( {{u}^{2}}\left( x \right) \right)\le 0$, i.e. $g\left( {{u}^{2}}\left( x \right) \right)\le \left( 1+\alpha  \right){{u}^{2}}\left( x \right)$ for $x\in {{\mathbb{R}}^{2}}$ if ${{u}^{2}}\left( x \right)\le {{s}_{\alpha }}$ for $x\in {{\mathbb{R}}^{2}}$, i.e. $\left\| u \right\|_{\infty }^{2}\le {{s}_{\alpha }}$. This completes the proof of Theorem~\ref{theomII}~(iii).

The rest of the proof of Theorem~\ref{theomII} is to show Theorem~\ref{theomII}~(ii), as follows. Suppose $-1<\alpha <-\frac{1}{2}$. Let ${{h}_{0}}\left( s \right)=h\left( s \right)-{{\rho }_{0}}s$ for $s>0$, where ${{\rho }_{0}}$ is a positive constant to be determined later. Then ${{h}_{0}}\left( 0 \right)=0$ and ${{h}_{0}}^{\prime }\left( s \right)=\frac{-1}{2\sqrt{1+s}}\left( 1+2\alpha +\frac{1}{1+s} \right)-{{\rho }_{0}}$ for $s>0$. To make ${{h}_{0}}^{\prime }\left( s \right)\le 0$ for $s>0$, we choose ${{\rho }_{0}}={{\left( -\frac{1+2\alpha }{3} \right)}^{\frac{3}{2}}}>0$ such that $-\left( 1+2\alpha  \right)-\frac{1}{1+s}\le 2{{\rho }_{0}}\sqrt{1+s}$ for $s>0$. To see this, we set $\omega \left( \tau  \right)=2{{\rho }_{0}}{{\tau }^{3}}+\left( 1+2\alpha  \right){{\tau }^{2}}+1$ for $\tau =\sqrt{1+s}>1$. We want to show $\omega \left( \tau  \right)\ge 0$ for $\tau >1$. By direct calculation, ${\omega }'\left( \tau  \right)=2\tau \left[ 3{{\rho }_{0}}\tau +\left( 1+2\alpha  \right) \right]$ and ${\omega }''\left( \tau  \right)=2\left[ 6{{\rho }_{0}}\tau +\left( 1+2\alpha  \right) \right]$, yielding ${\omega }'\left( -\frac{1+2\alpha }{3{{\rho }_{0}}} \right)=0$ and  ${\omega }''\left( -\frac{1+2\alpha }{3{{\rho }_{0}}} \right)=-\left( 1+2\alpha  \right)>0$. Moreover, $\omega \left( -\frac{1+2\alpha }{3{{\rho }_{0}}} \right)=\frac{1}{27}{{\left( 1+2\alpha  \right)}^{3}}\rho _{0}^{-2}+1\ge 0$ if $\rho _{0}^{2}\ge {{\left( -\frac{1+2\alpha }{3} \right)}^{3}}>0$. Note that the minimum of the function $\omega $ may be attained at $-\frac{1+2\alpha }{3{{\rho }_{0}}}={{\left( -\frac{1+2\alpha }{3} \right)}^{-\frac{1}{2}}}>1$ and $-\frac{1}{3}<\frac{1+2\alpha }{3}<0$, because $-1<\alpha <-\frac{1}{2}$. In particular, we set ${{\rho }_{0}}={{\left( -\frac{1+2\alpha }{3} \right)}^{\frac{3}{2}}}>0$. Then $\omega \left( \tau  \right)\ge 0$ for $\tau >1$, i.e., ${{h}_{0}}^{\prime }\left( s \right)\le 0$ for $s>0$ holds true. Therefore, ${{h}_{0}}\left( s \right)\le {{h}_{0}}\left( 0 \right)=0$ for $s>0$ and we have completed the proof of Theorem~\ref{theomII}~(ii).

\section{Proof of Theorem~\ref{theomI} }\label{saue1}
We multiply (\ref{eq-1.1}) by $u$ and integrate over ${{\mathbb{R}}^{2}}$. Then using integration by parts, we get
\BE\label{1.5}
-\int_{{{\mathbb{R}}^{2}}}{{{\left| \nabla u \right|}^{2}}}+\Gamma \int_{{{\mathbb{R}}^{2}}}{\left( 1-\frac{1}{1+{{u}^{2}}} \right)}{{u}^{2}}=\hat{\lambda }\left\| u \right\|_{2}^{2}\,,
\EE
where $\left\| u \right\|_{2}^{2}=\int_{{{\mathbb{R}}^{2}}}{{{u}^{2}}}$. On the other hand, we can also multiply (\ref{eq-1.1}) by $x\cdot \nabla u$ and integrate it over ${{\mathbb{R}}^{2}}$, where $x\cdot \nabla u=\sum\limits_{j=1}^{2}{{{x}_{j}}{{\partial }_{j}}u}$. Then using integration by parts, we derive the Pohozaev identity as follows
\BE\label{1.6}
\hat{\lambda }\left\| u \right\|_{2}^{2}=\Gamma \int_{{{\mathbb{R}}^{2}}}{\left[ {{u}^{2}}-\ln \left( 1+{{u}^{2}} \right) \right]}\,.
\EE
Combining (\ref{1.5}) and (\ref{1.6}), we have
$$
-\int_{{{\mathbb{R}}^{2}}}{{{\left| \nabla u \right|}^{2}}}+\Gamma \int_{{{\mathbb{R}}^{2}}}{\left( 1-\frac{1}{1+{{u}^{2}}} \right)}{{u}^{2}}=\Gamma \int_{{{\mathbb{R}}^{2}}}{\left[ {{u}^{2}}-\ln \left( 1+{{u}^{2}} \right) \right]}\,,
$$
implying
\BE\label{1.7}
\int_{{{\mathbb{R}}^{2}}}{{{\left| \nabla u \right|}^{2}}}=\Gamma \int_{{{\mathbb{R}}^{2}}}{\left( 1-\frac{1}{1+{{u}^{2}}} \right)}{{u}^{2}}-\Gamma \int_{{{\mathbb{R}}^{2}}}{\left[ {{u}^{2}}-\ln \left( 1+{{u}^{2}} \right) \right]}\,.
\EE
Hence
\BE\label{1.8}
\begin{array}{ll}
\beta &=\frac{\int_{{{\mathbb{R}}^{2}}}{{{\left| \nabla u \right|}^{2}}}}{-\Gamma \int_{{{\mathbb{R}}^{2}}}{\left[ {{u}^{2}}-\ln \left( 1+{{u}^{2}} \right) \right]}} \\
 & =1-\frac{\int_{{{\mathbb{R}}^{2}}}{\left( 1-\frac{1}{1+{{u}^{2}}} \right){{u}^{2}}}}{\int_{{{\mathbb{R}}^{2}}}{\left[ {{u}^{2}}-\ln \left( 1+{{u}^{2}} \right) \right]}} \\
\end{array}
\EE
It is obvious that $\beta <0$ because $u$ is nontrivial. To prove $\beta >-1$, we define a function $f\left( s \right)=\frac{\left( 1-\frac{1}{1+s} \right)s}{s-\ln \left( 1+s \right)}$ for $s>0$. Then we can show that ${f}'\left( s \right)<0$ for $s>0$, $\underset{s\to 0+}{\mathop{\lim }}\,f\left( s \right)=2$ and $\underset{s\to \infty }{\mathop{\lim }}\,f\left( s \right)=1$ (see Proposition~A.1 in Appendix~I). Consequently, $f\left( s \right)<2$ for $s>0$ and
\[\int_{{{\mathbb{R}}^{2}}}{\left( 1-\frac{1}{1+{{u}^{2}}} \right)}{{u}^{2}}=\int_{{{\mathbb{R}}^{2}}}{\left[ {{u}^{2}}-\ln \left( 1+{{u}^{2}} \right) \right]}f\left( {{u}^{2}} \right)\le 2\int_{{{\mathbb{R}}^{2}}}{\left[ {{u}^{2}}-\ln \left( 1+{{u}^{2}} \right) \right]}\]
which gives $\beta \ge -1$. Here we have used (\ref{1.8}).

Now we prove $\beta >-1$ by contradiction. Suppose $\beta =-1$, i.e.,
\BE\label{1.9}
\int_{{{\mathbb{R}}^{2}}}{\left( 1-\frac{1}{1+{{u}^{2}}} \right)}{{u}^{2}}=2\int_{{{\mathbb{R}}^{2}}}{\left[ {{u}^{2}}-\ln \left( 1+{{u}^{2}} \right) \right]}\,.
\EE
Due to the fact that $u\not\equiv 0$, the unique continuation theorem of equation (\ref{eq-1.1}) implies that there exists a ball ${{B}_{0}}={{B}_{{{R}_{0}}}}\left( {{x}_{0}} \right)$ with radius ${{R}_{0}}>0$ and center ${{x}_{0}}\in {{\mathbb{R}}^{2}}$ such that ${{u}^{2}}\left( x \right)>0$ for $x\in {{B}_{0}}$. Then there exist ${{\varepsilon }_{0}}>0$ and a smaller ball ${{B}_{1}}={{B}_{{{R}_{1}}}}\left( {{x}_{0}} \right)\subset \subset {{B}_{0}}$ such that ${{u}^{2}}\left( x \right)\ge {{\varepsilon }_{0}}>0$ for $x\in {{B}_{1}}$. Because of ${f}'\left( s \right)<0$ for $s>0$, $\underset{s\to 0+}{\mathop{\lim }}\,f\left( s \right)=2$ and $\underset{s\to \infty }{\mathop{\lim }}\,f\left( s \right)=1$ (see Proposition~A.1 in Appendix~I), there exists a ${{\delta }_{0}}>0$ such that $f\left( {{u}^{2}}\left( x \right) \right)\le 2-{{\delta }_{0}}$ for $x\in {{B}_{1}}$. Hence (\ref{1.9}) causes the following contradiction:
\begin{eqnarray*}
2\int_{{{\mathbb{R}}^{2}}}{\left[ {{u}^{2}}-\ln \left( 1+{{u}^{2}} \right) \right]} &=& \int_{{{\mathbb{R}}^{2}}}{\left( 1-\frac{1}{1+{{u}^{2}}} \right)}{{u}^{2}} \\
 &=& \int_{{{B}_{1}}}{\left( 1-\frac{1}{1+{{u}^{2}}} \right)}{{u}^{2}}+\int_{B_{1}^{c}}{\left( 1-\frac{1}{1+{{u}^{2}}} \right)}{{u}^{2}} \\
 &=& \int_{{{B}_{1}}}{f\left( {{u}^{2}} \right)}\left[ {{u}^{2}}-\ln \left( 1+{{u}^{2}} \right) \right]+\int_{B_{1}^{c}}{f\left( {{u}^{2}} \right)}\left[ {{u}^{2}}-\ln \left( 1+{{u}^{2}} \right) \right] \\
&\le& \left( 2-{{\delta }_{0}} \right)\int_{{{B}_{1}}}{\left[ {{u}^{2}}-\ln \left( 1+{{u}^{2}} \right) \right]}+2\int_{B_{1}^{c}}{\left[ {{u}^{2}}-\ln \left( 1+{{u}^{2}} \right) \right]} \\
&=& 2\int_{{{\mathbb{R}}^{2}}}{\left[ {{u}^{2}}-\ln \left( 1+{{u}^{2}} \right) \right]}-{{\delta }_{0}}\int_{{{B}_{1}}}{\left[ {{u}^{2}}-\ln \left( 1+{{u}^{2}} \right) \right]} \\
&<& 2\int_{{{\mathbb{R}}^{2}}}{\left[ {{u}^{2}}-\ln \left( 1+{{u}^{2}} \right) \right]}\,, \\
\end{eqnarray*}
where $B_{1}^{c}={{\mathbb{R}}^{2}}-{{B}_{1}}$ is the complement of ${{B}_{1}}$. Therefore, we have completed the proof of $-1<\beta <0$, i.e., (\ref{1.3}).

To prove (\ref{1.4}), we substitute (\ref{ba-1}) into (\ref{1.5}) and get
\BE\label{1.10}
\hat{\lambda }\left\| u \right\|_{2}^{2}=\beta \Gamma \int_{{{\mathbb{R}}^{2}}}{\left[ {{u}^{2}}-\ln \left( 1+{{u}^{2}} \right) \right]}+\Gamma \int_{{{\mathbb{R}}^{2}}}{\left( 1-\frac{1}{1+{{u}^{2}}} \right)}{{u}^{2}}=\Gamma \int_{{{\mathbb{R}}^{2}}}{g\left( {{u}^{2}} \right)}\,,
\EE
where $g\left( s \right)=\beta \left[ s-\ln \left( 1+s \right) \right]+\frac{{{s}^{2}}}{1+s}$  for $s>0$. Now we claim $g\left( s \right)<{{\left( 1+\frac{1}{2}\beta  \right)}^{2}}s$ for $s>0$. Let $h\left( s \right)=g\left( s \right)-{{\left( 1+\frac{1}{2}\beta  \right)}^{2}}s$ for $s>0$. Then by direct calculation, we get $h\left( s \right)=-\beta \ln \left( 1+s \right)+\frac{1}{1+s}-1-\frac{1}{4}{{\beta }^{2}}s$ and
\begin{eqnarray*}
{h}'\left( s \right) &=& \frac{-\beta }{1+s}-\frac{1}{{{\left( 1+s \right)}^{2}}}-\frac{1}{4}{{\beta }^{2}} \\
&=& -\frac{1}{{{\left( 1+s \right)}^{2}}}\left[ \frac{1}{4}{{\beta }^{2}}{{\left( 1+s \right)}^{2}}+\beta \left( 1+s \right)+1 \right] \\
&=& -\frac{{{\left[ \frac{1}{2}\beta \left( 1+s \right)+1 \right]}^{2}}}{{{\left( 1+s \right)}^{2}}}<0\,,
\end{eqnarray*}
for $s>0$. Therefore, $h\left( s \right)<h\left( 0 \right)=0$ i.e. $g\left( s \right)<{{\left( 1+\frac{1}{2}\beta  \right)}^{2}}s$ for $s>0$, which can be put into (\ref{1.10}) to get (\ref{1.4}) and, thus, complete the proof of Theorem~\ref{theomI}.

\section{Proof of Theorem~\ref{theomIII}}\label{saue2}
We first prove the upper bound estimate of the ground state energy ${{e}_{\Gamma }}$.
\BL\label{saue3.1} Given the same assumptions as in Theorem~\ref{theomIII}, we have
${{e}_{\Gamma }}\le -\frac{\Gamma }{2}\left( 1+{{o}_{\Gamma }}\left( 1 \right) \right)$ as $\Gamma \to\infty $., where ${{o}_{\Gamma }}\left( 1 \right)$  is a small quantity tending to zero as $\Gamma $ goes to infinity.
\EL
\begin{proof}
Let ${{u}_{\tau }}\left( x \right)=\frac{1}{\tau }U\left( \frac{x}{\tau } \right)$  for $x\in {{\mathbb{R}}^{2}}$, $\tau>0$, where  $\text{U}\in {{\text{H}}^{1}}\left( {{\mathbb{R}}^{2}} \right)$, ${{\left\| U \right\|}_{2}}=1$ and $U>0$ in  ${{\mathbb{R}}^{2}}$. Then
${{\left\| {{u}_{\tau }} \right\|}_{2}}={{\left\| U \right\|}_{2}}=1$, $\int_{{{\mathbb{R}}^{2}}}{{{\left| \nabla {{u}_{\tau }}\left( x \right) \right|}^{2}}dx}={{\tau }^{-2}}{{\int_{{{\mathbb{R}}^{2}}}{\left| \nabla U\left( x \right) \right|}}^{2}}dx$ and
\begin{eqnarray*}
E\left[ {{u}_{\tau }} \right]&=&\frac{1}{2}\int_{{{\mathbb{R}}^{2}}}{{{\left| \nabla {{u}_{\tau }} \right|}^{2}}-\Gamma \left[ u_{\tau }^{2}-\ln \left( 1+u_{\tau }^{2} \right) \right]dx} \\
&=&\frac{1}{2}\int_{{{\mathbb{R}}^{2}}}{{{\left| \nabla {{u}_{\tau }} \right|}^{2}}dx}-\frac{\Gamma }{2}\int_{{{\mathbb{R}}^{2}}}{u_{\tau }^{2}dx}+\Gamma \int_{{{\mathbb{R}}^{2}}}{\ln \left( 1+u_{\tau }^{2} \right)dx} \\
&=&\frac{1}{2}{{\tau }^{-2}}\int_{{{\mathbb{R}}^{2}}}{{{\left| \nabla U \right|}^{2}}dx}-\frac{\Gamma }{2}+\Gamma \int_{{{\mathbb{R}}^{2}}}{{{\tau }^{2}}\ln \left( 1+{{\tau }^{-2}}{{U}^{2}}\left( y \right) \right)dy} \\
&=& -\frac{\Gamma }{2}\left( 1+{{o}_{\Gamma }}\left( 1 \right) \right) \,, \\
\end{eqnarray*}
as $\tau \sim {{\left( \ln \Gamma  \right)}^{-1/2}}$ and $\Gamma \to \infty $, where $y=\frac{x}{\tau }$ and ${{o}_{\Gamma }}\left( 1 \right)$ is a small quantity tending to zero as $\Gamma $ goes to infinity. Here we have used the fact ${{\tau }^{2}}\ln \left( 1+{{\tau }^{-2}}{{U}^{2}} \right)\le C{{U}^{2}}\in {{L}^{1}}\left( {{\mathbb{R}}^{2}} \right)$ for some constant $C>0$ (independent of $U$ and $\tau $) and $0<{{\tau }^{2}}\ln \left( 1+{{\tau }^{-2}}{{U}^{2}}\left( y \right) \right)\le {{\tau }^{2}}\ln \left( 1+{{\tau }^{-2}}\left\| U \right\|_{\infty }^{2} \right)\to 0$ as $\tau \to 0$ for $y\in {{\mathbb{R}}^{2}}$. Hence, by the Dominated Convergence Theorem,
$$
\int_{{{\mathbb{R}}^{2}}}{{{\tau }^{2}}\ln \left( 1+{{\tau }^{-2}}{{U}^{2}}\left( y \right) \right)dy}={{o}_{\Gamma }}\left( 1 \right)
\quad\hbox{ as }\quad \Gamma\to\infty\,.
$$
Note that $\tau \sim {{\left( \ln \Gamma  \right)}^{-1/2}}$ as $\Gamma \to \infty $. Therefore, ${{e}_{\Gamma }}\le E\left[ {{u}_{\tau }} \right]=-\frac{\Gamma }{2}\left( 1+o\left( 1 \right) \right)$ as $\Gamma \to \infty$ and we have completed the proof of Lemma~\ref{saue3.1}.
\end{proof}

For the lower bound estimate of ${{e}_{\Gamma }}$, it is obvious that
\begin{eqnarray*}
{{e}_{\Gamma }}&=& E\left[ {{u}_{\Gamma }} \right]=\frac{1}{2}\int_{{{\mathbb{R}}^{2}}}{{{\left| \nabla {{u}_{\Gamma }} \right|}^{2}}-\Gamma \left[ u_{\Gamma }^{2}-\ln \left( 1+u_{\Gamma }^{2} \right) \right]dx} \\
 &=&\frac{1}{2}\int_{{{\mathbb{R}}^{2}}}{{{\left| \nabla {{u}_{\Gamma }} \right|}^{2}}dx}-\frac{\Gamma }{2}\int_{{{\mathbb{R}}^{2}}}{u_{\Gamma }^{2}dx}+\Gamma \int_{{{\mathbb{R}}^{2}}}{\ln \left( 1+u_{\Gamma }^{2} \right)dx} \\
 &\ge& -\frac{\Gamma }{2}\int_{{{\mathbb{R}}^{2}}}{u_{\Gamma }^{2}dx}=-\frac{\Gamma }{2}\,, \\
\end{eqnarray*}
where ${{u}_{\Gamma }}$ is the ground state (energy minimizer) of ${{e}_{\Gamma }}$ under the ${{L}^{2}}$-norm constraint ${{\left\| {{u}_{\Gamma }} \right\|}_{2}}=1$. Consequently, by Lemma~\ref{saue3.1}, we obtain ${{e}_{\Gamma }}=-\frac{\Gamma }{2}\left( 1+{{o}_{\Gamma }}\left( 1 \right) \right)$ as $\Gamma \to \infty$. To get a further estimate of ${{e}_{\Gamma }}$, we need the following lemmas.

\BL\label{saue3.2}
Continuing from Lemma~\ref{saue3.1}, we have that the ratio ${{\beta }_{\Gamma }}=\frac{\int_{{{\mathbb{R}}^{2}}}{{{\left| \nabla {{u}_{\Gamma }} \right|}^{2}}dx}}{-\Gamma \int_{{{\mathbb{R}}^{2}}}{\left[ u_{\Gamma }^{2}-\ln \left( 1+u_{\Gamma }^{2} \right) \right]dx}}\to 0$  as $\Gamma \to \infty $, where ${{u}_{\Gamma }}$ is the ground state (energy minimizer) of ${{e}_{\Gamma }}$.
\EL
\begin{proof}
We prove by contradiction. Suppose that ${{\beta }_{\Gamma }}$ does not approach zero as $\Gamma $ goes to infinity. Then by (\ref{1.3}) of Theorem~\ref{theomI}, we may assume ${{\beta }_{\Gamma }}\to -{{c}_{0}}$  as $\Gamma \to \infty $, where $0<{{c}_{0}}\le 1$ is a constant. Hence
\begin{align}\label{saue2-eq1}
& {{e}_{\Gamma }}=\frac{1}{2}\int_{{{\mathbb{R}}^{2}}}{{{\left| \nabla {{u}_{\Gamma }} \right|}^{2}}} -\frac{\Gamma }{2}\int_{{{\mathbb{R}}^{2}}}{\left[ u_{\Gamma }^{2}-\ln \left( 1+u_{\Gamma }^{2} \right) \right]} \\
& =-\frac{\Gamma }{2}\left( {{\beta }_{\Gamma }}+1 \right)\int_{{{\mathbb{R}}^{2}}}{\left[ u_{\Gamma }^{2}-\ln \left( 1+u_{\Gamma }^{2} \right) \right]}  \nonumber \\
& \ge -\frac{\Gamma }{2}\left( {{\beta }_{\Gamma }}+1 \right)\int_{{{\mathbb{R}}^{2}}}{u_{\Gamma }^{2}}  \nonumber \\
& \text{=}-\frac{\Gamma }{2}\left( {{\beta }_{\Gamma }}+1 \right)\,. \nonumber \\
\end{align}
Here we have used the fact that $-1<{{\beta }_{\Gamma }}<0$  and ${{\left\| {{u}_{\Gamma }} \right\|}_{2}}=1$. Combining Lemma~\ref{saue3.1} and (\ref{saue2-eq1}), we have $-\frac{1}{2}\left( 1+{{o}_{\Gamma }}\left( 1 \right) \right)\ge -\frac{1}{2}\left( {{\beta }_{\Gamma }}+1 \right)$, i.e. ${{\beta }_{\Gamma }}\ge 0$ as $\Gamma \to \infty $, which contradicts with ${{\beta }_{\Gamma }}\to -{{c}_{0}}\in \left[ -1,0 \right)$ as $\Gamma \to \infty $.
Therefore, we have completed the proof of Lemma~\ref{saue3.2}.
\end{proof}

\BL\label{saue2.0}
Continuing from Lemma~\ref{saue3.1}, we have that
${{e}_{\Gamma }}$ is decreasing to $\Gamma $ for $\Gamma >{{T}_{1}}$.
\EL
\begin{proof}
 Let ${{u}_{\Gamma }}$ be the energy minimizer (ground state) of ${{e}_{\Gamma }}$  for $\Gamma >{{T}_{2}}$. Then
\begin{eqnarray*}
 2{{e}_{{{\Gamma }_{1}}}}&=&{{\int_{{{\mathbb{R}}^{2}}}{\left| \nabla {{u}_{{{\Gamma }_{1}}}} \right|}^{2}}}-{{\Gamma }_{1}}\int_{{{\mathbb{R}}^{2}}}{\left[ u_{{{\Gamma }_{1}}}^{2}-\ln \left( 1+u_{{{\Gamma }_{1}}}^{2} \right) \right]} \\
 &\ge& {{\int_{{{\mathbb{R}}^{2}}}{\left| \nabla {{u}_{{{\Gamma }_{1}}}} \right|}^{2}}}-{{\Gamma }_{2}}\int_{{{\mathbb{R}}^{2}}}{\left[ u_{{{\Gamma }_{1}}}^{2}-\ln \left( 1+u_{{{\Gamma }_{1}}}^{2} \right) \right]} \\
 &\ge& 2{{e}_{{{\Gamma }_{2}}}}\,, \\
\end{eqnarray*}
for ${{\Gamma }_{2}}>{{\Gamma }_{1}}>{{T}_{1}}>0$. Hence ${{e}_{\Gamma }}$ is decreasing to $\Gamma $ and we have completed the proof of Lemma~\ref{saue2.0}.
\end{proof}

\BL\label{saue2.1}
Continuing from Lemma~\ref{saue3.1}, we have that
${{e}_{\Gamma +1}}-{{e}_{\Gamma }}\ge -\frac{1}{2}\frac{{{{\hat{\lambda }}}_{\Gamma \text{+}1}}}{\Gamma +1}$ for $\Gamma >{{T}_{1}}$ , where ${{u}_{\Gamma +1}}$ is the energy minimizer (ground state) of ${{e}_{\Gamma +1}}$ and ${{\hat{\lambda }}_{\Gamma +1}}$ is the associated eigenvalue of ${{u}_{\Gamma +1}}$.
\EL
\begin{proof}
It is obvious that for $\Gamma >{{T}_{1}}$ ,
\begin{eqnarray*}
2{{e}_{\Gamma +1}}&=&\int_{{{\mathbb{R}}^{2}}}{{{\left| \nabla {{u}_{\Gamma +1}} \right|}^{2}}-\left( \Gamma +1 \right)\int_{{{\mathbb{R}}^{2}}}{\left[ u_{\Gamma +1}^{2}-\ln \left( 1+u_{\Gamma +1}^{2} \right) \right]}} \\
 &=& \int_{{{\mathbb{R}}^{2}}}{{{\left| \nabla {{u}_{\Gamma +1}} \right|}^{2}}-\Gamma }\int_{{{\mathbb{R}}^{2}}}{\left[ u_{\Gamma +1}^{2}-\ln \left( 1+u_{\Gamma +1}^{2} \right) \right]}-\int_{{{\mathbb{R}}^{2}}}{\left[ u_{\Gamma +1}^{2}-\ln \left( 1+u_{\Gamma +1}^{2} \right) \right]} \\
 &\ge& 2{{e}_{\Gamma }}-\int_{{{\mathbb{R}}^{2}}}{\left[ u_{\Gamma +1}^{2}-\ln \left( 1+u_{\Gamma +1}^{2} \right) \right]}\,. \\
\end{eqnarray*}
Therefore, we may use the Pohozaev identity (\ref {1.6}) to complete the proof of Lemma~\ref{saue2.1} just as before.
\end{proof}

Now we want to prove (\ref{lb-est1}), i.e., ${{e}_{\Gamma }}\ge -\frac{\Gamma }{2}+\sigma \frac{{{T}_{1}}}{2}\ln \Gamma +{{C}_{0}}$ for $\sigma \in \left( 0,1 \right)$ and $\Gamma >{{\Gamma }_{\sigma }}$ sufficiently large, where ${{\Gamma }_{\sigma }}$ is a positive constant depending on $\sigma $ and ${{C}_{0}}$ is a constant independent of $\Gamma $. By Lemma~\ref{saue2.1} and (\ref{1.4}) of Theorem~\ref{theomI}, we have
\BE\label{saue-eq3.1}
{{e}_{\Gamma +1}}-{{e}_{\Gamma }}\ge -\frac{1}{2}{{\left( 1+\frac{1}{2}{{\beta }_{\Gamma +1}} \right)^{2}}}=-\frac{1}{2}\left[ 1+{{\beta }_{\Gamma +1}}\left( 1+\frac{1}{4}{{\beta }_{\Gamma +1}} \right) \right]\,,
\quad\hbox{ for }\quad \Gamma >{{T}_{1}}\,.
\EE
From Lemma~\ref{saue3.2}, $0>{{\beta }_{\Gamma }}\to 0$ as $\Gamma \to \infty $, implying that for all $\sigma \in \left( 0,1 \right)$, there exists a positive constant ${{\Gamma }_{\sigma }}$ sufficiently large such that $1+\frac{1}{4}{{\beta }_{\Gamma +1}}>\sigma $ for $\Gamma >{{\Gamma }_{\sigma }}$. Hence, (\ref{saue-eq3.1}) becomes
\BE\label{saue-eq3.2}
{{e}_{\Gamma +1}}-{{e}_{\Gamma }}\ge -\frac{1}{2}\left( 1+\sigma {{\beta }_{\Gamma +1}} \right)\ge -\frac{1}{2}\left( 1-\sigma \frac{{{T}_{1}}}{\Gamma +1} \right) \quad\hbox{ for }\quad \sigma \in \left( 0,1 \right)\:\hbox{ and }\: \Gamma >{{\Gamma }_{\sigma }}\,.
\EE
Here we have used the fact that ${{\beta }_{\Gamma +1}}\le -\frac{{{T}_{1}}}{\Gamma +1}$ due to ${{\beta }_{s}}=-\frac{1}{s}\frac{\int_{{{\mathbb{R}}^{2}}}{{{\left| \nabla {{u}_{s}} \right|}^{2}}}}{\int_{{{\mathbb{R}}^{2}}}{\left[ u_{s}^{2}-\ln \left( 1+u_{s}^{2} \right) \right]}}\le -\frac{{{T}_{1}}}{s}$ for $s>{{T}_{1}}$ and due to the definition of ${{T}_{1}}$.

Fix $\sigma \in \left( 0,1 \right)$ arbitrarily and let $N\in \mathbb{N}$ and $N>{{\Gamma }_{\sigma }}$. Then (\ref{saue-eq3.2}) gives
$$
{{e}_{k+1}}-{{e}_{k}}\ge -\frac{1}{2}\left( 1-\sigma \frac{{{T}_{1}}}{k+1} \right)\quad\hbox{ for }\quad k=N,N+1,N+2,\cdots\,.
$$
Hence, for $n\in \mathbb{N}$,
\begin{eqnarray*}
{{e}_{N+n}}-{{e}_{N}}&=&\sum\limits_{j=0}^{n-1}{\left( {{e}_{N+j+1}}-{{e}_{N+j}} \right)} \\
&\ge& -\frac{1}{2}\sum\limits_{j=0}^{n-1}{\left( 1-\sigma \frac{{{T}_{1}}}{N+j+1} \right)} \\
&=& -\frac{n}{2}+\sigma \frac{{{T}_{1}}}{2}\sum\limits_{k=1}^{n}{\frac{1}{N+k}} \\
&\ge& -\frac{n}{2}+\sigma \frac{{{T}_{1}}}{2}\int_{N+1}^{N+n}{\frac{1}{t}dt} \\
&=&-\frac{n}{2}+\sigma \frac{{{T}_{1}}}{2}\left[ \ln \left( N+n \right)-\ln \left( N+1 \right) \right] \,,
\end{eqnarray*}
yielding ${{e}_{\left[ \Gamma  \right]+1}}\ge -\frac{1}{2}\left[ \Gamma  \right]+\sigma \frac{{{T}_{1}}}{2}\ln \left( \left[ \Gamma  \right]+1 \right)+{{C}_{N}}$  for $\Gamma >N$ sufficiently large, where $\left[ \Gamma  \right]=\sup \left\{ k\in \mathbb{N}:k\le \Gamma  \right\}$ and where we set $\left[ \Gamma  \right]+1=N+n$. Consequently, we get ${{e}_{\Gamma }}\ge {{e}_{\left[ \Gamma  \right]+1}}\ge -\frac{1}{2}\Gamma +\sigma \frac{{{T}_{1}}}{2}\ln \Gamma +{{C}_{0}}$ and have completed the proof of (\ref{lb-est1}) because $\left[ \Gamma  \right]\le \Gamma \le \left[ \Gamma  \right]+1$ and ${{e}_{\Gamma }}$ is decreasing to $\Gamma $ for $\Gamma >{{T}_{1}}$. Here ${{C}_{0}}$ is a constant independent of $\Gamma $. Therefore, we have completed\emph{} the proof of Theorem~\ref{theomIII}.

\section{Proof of Theorem~\ref{theomIV}}\label{sqrt2}

We first prove the upper bound estimate of the ground state energy ${{\tilde{e}}_{\Gamma }}$.
\BL\label{sqrt-clm1} Under the same assumptions of Theorem~\ref{theomIV}, we have
${{\tilde{e}}_{\Gamma }}\le -\frac{\Gamma }{2}\left( 1+o_\Gamma\left( 1 \right) \right)$ as $\Gamma \to \infty $, where $o_\Gamma\left( 1 \right)$ is a small quantity tending to zero as $\Gamma $ goes to infinity.
\EL
\begin{proof}
Let ${{u}_{\tau }}\left( x \right)=\frac{1}{\tau }U\left( \frac{x}{\tau } \right)$ for $x\in {{\mathbb{R}}^{2}}$, $\tau>0$, where $\text{U}\in {{{H}}^{1}}\left( {{\mathbb{R}}^{2}} \right)\bigcap {{L}^{1}}\left( {{\mathbb{R}}^{2}} \right)$, ${{\left\| U \right\|}_{2}}=1$ and $U>0$ in ${{\mathbb{R}}^{2}}$. Then ${{\left\| {{u}_{\tau }} \right\|}_{2}}={{\left\| U \right\|}_{2}}=1$, $\int_{{{\mathbb{R}}^{2}}}{{{\left| \nabla {{u}_{\tau }}\left( x \right) \right|}^{2}}dx}={{\tau }^{-2}}{{\int_{{{\mathbb{R}}^{2}}}{\left| \nabla U\left( x \right) \right|}}^{2}}dx$ and
\begin{eqnarray*}
E\left[ {{u}_{\tau }} \right] &=& \frac{1}{2}\int_{{{\mathbb{R}}^{2}}}{{{\left| \nabla {{u}_{\tau }} \right|}^{2}}-\Gamma {{\left( \sqrt{1+u_{\tau }^{2}}-1 \right)}^{2}}dx} \\
 &=& \frac{1}{2}\int_{{{\mathbb{R}}^{2}}}{{{\left| \nabla {{u}_{\tau }} \right|}^{2}}dx}-\frac{\Gamma }{2}\int_{{{\mathbb{R}}^{2}}}{u_{\tau }^{2}dx}-\Gamma \int_{{{\mathbb{R}}^{2}}}{\left( 1-\sqrt{1+u_{\tau }^{2}} \right)dx} \\
 &=& \frac{1}{2}{{\tau }^{-2}}\int_{{{\mathbb{R}}^{2}}}{{{\left| \nabla U \right|}^{2}}dx}-\frac{\Gamma }{2}+\Gamma \int_{{{\mathbb{R}}^{2}}}{\frac{u_{\tau }^{2}}{1+\sqrt{1+u_{\tau }^{2}}}dx} \\
 &=& \frac{1}{2}{{\tau }^{-2}}\int_{{{\mathbb{R}}^{2}}}{{{\left| \nabla U \right|}^{2}}dx}-\frac{\Gamma }{2}+\Gamma \tau \int_{{{\mathbb{R}}^{2}}}{\frac{{{U}^{2}}}{\tau +\sqrt{{{\tau }^{2}}+{{U}^{2}}}}dx} \\
 &=& -\frac{\Gamma }{2}\left( 1+o_\Gamma\left( 1 \right) \right) \\
\end{eqnarray*}
as $\tau \sim {{\left( \ln \Gamma  \right)}^{-1/2}}$ and $\Gamma \to \infty $, where $o_\Gamma\left( 1 \right)$ is a small quantity tending to zero as $\Gamma $ goes to infinity. Here we have used the fact that $\tau \int_{{{\mathbb{R}}^{2}}}{\frac{{{U}^{2}}}{\tau +\sqrt{{{\tau }^{2}}+{{U}^{2}}}}dx}\le \tau \int_{{{\mathbb{R}}^{2}}}{Udx}\to 0$ because $U\in {{L}^{1}}\left( {{\mathbb{R}}^{2}} \right)$ and $\tau \sim {{\left( \ln \Gamma \right)}^{-1/2}}\to 0$ as $\Gamma \to \infty $. Therefore, we get ${{\tilde{e}}_{\Gamma }}\le E\left[ {{u}_{\tau }} \right]=-\frac{\Gamma }{2}\left( 1+o_\Gamma\left( 1 \right) \right)$ as $\Gamma \to \infty $ and have completed the proof of Lemma~\ref{sqrt-clm1}.
\end{proof}

For the lower bound estimate of ${{\tilde{e}}_{\Gamma }}$, it is obvious that
\begin{eqnarray*}
{{{\tilde{e}}}_{\Gamma }} &=& E\left[ {{u}_{\Gamma }} \right]=\frac{1}{2}\int_{{{\mathbb{R}}^{2}}}{{{\left| \nabla {{u}_{\Gamma }} \right|^{2}}}-\Gamma {{\left( \sqrt{1+u_{\Gamma }^{2}}-1 \right)}^{2}}dx} \\
 &=& \frac{1}{2}\int_{{{\mathbb{R}}^{2}}}{{{\left| \nabla {{u}_{\Gamma }} \right|}^{2}}dx}-\frac{\Gamma }{2}\int_{{{\mathbb{R}}^{2}}}{u_{\Gamma }^{2}dx}+\Gamma \int_{{{\mathbb{R}}^{2}}}{\left( \sqrt{1+u_{\Gamma }^{2}}-1 \right)dx} \\
 &\ge& -\frac{\Gamma }{2}\int_{{{\mathbb{R}}^{2}}}{u_{\Gamma }^{2}dx}=-\frac{\Gamma }{2}\,, \\
\end{eqnarray*}
where ${{u}_{\Gamma }}$ is the ground state (energy minimizer) of ${{\tilde{e}}_{\Gamma }}$ under the ${{L}^{2}}$-norm constraint ${{\left\| {{u}_{\Gamma }} \right\|}_{2}}=1$. Consequently, by Lemma~\ref{sqrt-clm1}, we obtain ${{\tilde{e}}_{\Gamma }}=-\frac{\Gamma }{2}\left( 1+{{o}_{\Gamma }}\left( 1 \right) \right)$ as $\Gamma \to \infty$. To get a further estimate of ${{\tilde{e}}_{\Gamma }}$, we need the following lemmas.

\BL\label{sqrt-clm2} Under the same assumptions of Theorem~\ref{theomIV}, we have that the ratio ${{\alpha }_{\Gamma }}=\frac{\int_{{{\mathbb{R}}^{2}}}{{{\left| \nabla {{u}_{\Gamma }} \right|^{2}}}dx}}{-\Gamma \int_{{{\mathbb{R}}^{2}}}{{{\left( \sqrt{1+u_{\Gamma }^{2}}-1 \right)^{2}}}dx}}\to 0$ as $\Gamma \to \infty $, where ${{u}_{\Gamma }}$ is the ground state (energy minimizer) of ${{\tilde{e}}_{\Gamma }}$.
\EL
\begin{proof}
We prove by contradiction. Suppose that ${{\alpha }_{\Gamma }}$ may not approach zero as $\Gamma $ goes to infinity. Then by (\ref {a.3}) of Theorem~\ref {theomII}, we may assume ${{\alpha }_{\Gamma }}\to -{{c}_{0}}$ as $\Gamma \to \infty $, where $0<{{c}_{0}}\le 1$ is a constant. Hence
\begin{eqnarray}\label{sqrt-eq1}
{{{\tilde{e}}}_{\Gamma }} &=& \frac{1}{2}\int_{{{\mathbb{R}}^{2}}}{{{\left| \nabla {{u}_{\Gamma }} \right|}^{2}}}-\frac{\Gamma }{2}\int_{{{\mathbb{R}}^{2}}}{{{\left( \sqrt{1+u_{\Gamma }^{2}}-1 \right)}^{2}}} \\
&=& -\frac{\Gamma }{2}\left( {{\alpha }_{\Gamma }}+1 \right)\int_{{{\mathbb{R}}^{2}}}{{{\left( \sqrt{1+u_{\Gamma }^{2}}-1 \right)}^{2}}}  \nonumber \\
&=& -\frac{\Gamma }{2}\left( {{\alpha }_{\Gamma }}+1 \right)\left[ \int_{{{\mathbb{R}}^{2}}}{u_{\Gamma }^{2}}+2\int_{{{\mathbb{R}}^{2}}}{\left( 1-\sqrt{1+u_{\Gamma }^{2}} \right)} \right]  \nonumber \\
&\ge& -\frac{\Gamma }{2}\left( {{\alpha }_{\Gamma }}+1 \right)\,. \nonumber
\end{eqnarray}
Here we have used the fact that $-1<{{\alpha }_{\Gamma }}<0$ and ${{\left\| {{u}_{\Gamma }} \right\|}_{2}}=1$. Combining Lemma~\ref{sqrt-clm1} and (\ref{sqrt-eq1}), we have $-\frac{1}{2}\left( 1+o\left( 1 \right) \right)\ge -\frac{1}{2}\left( {{\alpha }_{\Gamma }}+1 \right)$ i.e. ${{\alpha }_{\Gamma }}\ge 0$ as $\Gamma \to \infty $, contradicting ${{\alpha }_{\Gamma }}\to -{{c}_{0}}\in \left[ -1,0 \right)$ as $\Gamma \to \infty $.
\end{proof}

\BR\label{sqrt-rk1}
Lemma~\ref{sqrt-clm2} shows that ${{\alpha }_{\Gamma }}\to 0$ as $\Gamma \to \infty $ so the condition $-\frac{1}{2}\le {{\alpha }_{\Gamma }}<0$ of Theorem~\ref{theomII}~(i) can be satisfied as $\Gamma$ becomes sufficiently large. Consequently, we obtain the eigenvalue estimate
\BE\label{sqrt-eq2}
0<{{\tilde{\lambda }}_{\Gamma }}\le \Gamma \left( 1+{{\alpha }_{\Gamma }} \right) \hbox{ for }\quad \Gamma >{{\Gamma }_{0}}\,,
\EE
where ${{\Gamma }_{0}}$ is a positive constant and ${{\tilde{\lambda }}_{\Gamma }}$ is the eigenvalue of the ground state ${{u}_{\Gamma }}$ with ground state energy ${{\tilde{e}}_{\Gamma }}$.
\ER

\BL\label{sqrt-clm2.0}
Under the same assumptions as in Theorem~\ref{theomIV},
${{\tilde{e}}_{\Gamma }}$ is decreasing to $\Gamma $ for $\Gamma >{{T}_{2}}$.
\EL

\begin{proof}
 Let ${{u}_{\Gamma }}$ be the energy minimizer (ground state) of ${{\tilde{e}}_{\Gamma }}$ for $\Gamma >{{T}_{2}}$. Then
\begin{eqnarray*}
 {{{\tilde{e}}}_{{{\Gamma }_{1}}}} &=& \frac{1}{2}{{\int_{{{\mathbb{R}}^{2}}}{\left| \nabla {{u}_{{{\Gamma }_{1}}}} \right|}^{2}}}-{{\Gamma }_{1}}{{\left( \sqrt{1+u_{{{\Gamma }_{1}}}^{2}}-1 \right)^{2}}}dx \\
&\ge& \frac{1}{2}{{\int_{{{\mathbb{R}}^{2}}}{\left| \nabla {{u}_{{{\Gamma }_{1}}}} \right|}^{2}}}-{{\Gamma }_{2}}{{\left( \sqrt{1+u_{{{\Gamma }_{1}}}^{2}}-1 \right)^{2}}}dx \\
&\ge& {{{\tilde{e}}}_{{{\Gamma }_{2}}}}\,,
\end{eqnarray*}
for ${{\Gamma }_{2}}>{{\Gamma }_{1}}>{{T}_{2}}>0$. Hence, ${{\tilde{e}}_{\Gamma }}$ is decreasing with respect to $\Gamma$ and we have completed the proof of Lemma~\ref{sqrt-clm2.0}.
\end{proof}

\BL\label{sqrt-clm2.1}
Under the same assumptions as in Theorem~\ref{theomIV},
${{\tilde{e}}_{\Gamma +1}}-{{\tilde{e}}_{\Gamma }}\ge -\frac{1}{2}\frac{{{{\tilde{\lambda }}}_{\Gamma +1}}}{\Gamma +1}$ for $\Gamma >{{T}_{2}}$, where ${{u}_{\Gamma +1}}$ is the energy minimizer (ground state) of ${{\tilde{e}}_{\Gamma +1}}$.
\EL

\begin{proof}
It is obvious that for $\Gamma >{{T}_{2}}$,
\begin{eqnarray*}
{{{\tilde{e}}}_{\Gamma +1}} &=& \frac{1}{2}\int_{{{\mathbb{R}}^{2}}}{{{\left| \nabla {{u}_{\Gamma +1}} \right|^{2}}}-\left( \Gamma +1 \right){{\left( \sqrt{1+u_{\Gamma +1}^{2}}-1 \right)^{2}}}} \\
 &=& \frac{1}{2}\int_{{{\mathbb{R}}^{2}}}{{{\left| \nabla {{u}_{\Gamma +1}} \right|^{2}}}-\Gamma {{\left( \sqrt{1+u_{\Gamma +1}^{2}}-1 \right)^{2}}}}-\frac{1}{2}\int_{{{\mathbb{R}}^{2}}}{{{\left( \sqrt{1+u_{\Gamma +1}^{2}}-1 \right)^{2}}}} \\
 &\ge& {{{\tilde{e}}}_{\Gamma }}-\frac{1}{2}\int_{{{\mathbb{R}}^{2}}}{{{\left( \sqrt{1+u_{\Gamma +1}^{2}}-1 \right)^{2}}}}\,.
\end{eqnarray*}
Therefore, we can use the Pohozaev identity (\ref {a.7}) to complete the proof of Lemma~\ref{sqrt-clm2.1}.
\end{proof}

Now we are ready to prove (\ref{lb-est2}), i.e., ${{\tilde{e}}_{\Gamma }}\ge -\frac{\Gamma }{2}+\frac{{{T}_{2}}}{2}\ln \Gamma +{{C}_{0}}$ for $\Gamma >{{\Gamma }_{0}}$ sufficiently large, where ${{C}_{0}}$ is a constants independent of $\Gamma $. By Lemma~\ref{sqrt-clm2.1} and (\ref{sqrt-eq2}), we have
\BE\label{sqrt-eq3.1}
{{\tilde{e}}_{\Gamma +1}}-{{\tilde{e}}_{\Gamma }}\ge -\frac{1}{2}\left( 1+{{\alpha }_{\Gamma +1}} \right)\ge -\frac{1}{2}\left( 1-\frac{{{T}_{2}}}{\Gamma +1} \right)
\EE
for $\Gamma >{{\Gamma }_{0}}$. Here we have used the fact that ${{\alpha }_{\Gamma +1}}\le -\frac{{{T}_{2}}}{\Gamma +1}$ from
$$
{{\alpha }_{s}}=-\frac{1}{s}\frac{\int_{{{\mathbb{R}}^{2}}}{{{\left| \nabla {{u}_{s}} \right|}^{2}}}}{\int_{{{\mathbb{R}}^{2}}}{{{\left( \sqrt{1+u_{s}^{2}}-1 \right)}^{2}}}}\le -\frac{{{T}_{2}}}{s}
\quad{ for }\quad s>{{T}_{2}}\,,
$$
due to the definition of ${{T}_{2}}$. Fix $N\in \mathbb{N}$ and $N>{{\Gamma }_{0}}$. Then (\ref{sqrt-eq3.1}) gives
\[{{\tilde{e}}_{k+1}}-{{\tilde{e}}_{k}}\ge -\frac{1}{2}\left( 1-\frac{{{T}_{2}}}{k+1} \right)  \quad\hbox{ for }\quad k=N,N+1,N+2,\cdots  \]
Consequently, for $n\in \mathbb{N}$,
\begin{eqnarray*}
{{{\tilde{e}}}_{N+n}}-{{{\tilde{e}}}_{N}} &=& \sum\limits_{j=0}^{n-1}{\left( {{{\tilde{e}}}_{N+j+1}}-{{{\tilde{e}}}_{N+j}} \right)} \\
 &\ge& -\frac{1}{2}\sum\limits_{j=0}^{n-1}{\left( 1-\frac{{{T}_{2}}}{N+j+1} \right)} \\
 &=& -\frac{n}{2}+\frac{{{T}_{2}}}{2}\sum\limits_{k=1}^{n}{\frac{1}{N+k}} \\
 &\ge& -\frac{n}{2}+\frac{{{T}_{2}}}{2}\int_{N+1}^{N+n}{\frac{1}{t}dt} \\
 &=& -\frac{n}{2}+\frac{{{T}_{2}}}{2}\left[ \ln \left( N+n \right)-\ln \left( N+1 \right) \right]\,,
\end{eqnarray*}
yielding ${{\tilde{e}}_{\left[ \Gamma  \right]+1}}\ge -\frac{1}{2}\left[ \Gamma  \right]+\frac{{{T}_{2}}}{2}\ln \left( \left[ \Gamma  \right]+1 \right)+{{C}_{N}}$ for $\Gamma >N$ sufficiently large, where $\left[ \Gamma  \right]=\sup \left\{ k\in \mathbb{N}:k\le \Gamma  \right\}$ and where we set $\left[ \Gamma  \right]+1=N+n$. Thus we have ${{\tilde{e}}_{\Gamma }}\ge {{\tilde{e}}_{\left[ \Gamma  \right]+1}}\ge -\frac{1}{2}\Gamma +\frac{{{T}_{2}}}{2}\ln \Gamma +{{C}_{0}}$ because $\left[ \Gamma  \right]\le \Gamma \le \left[ \Gamma  \right]+1$ and ${{\tilde{e}}_{\Gamma }}$ is decreasing to $\Gamma $ for $\Gamma >{{T}_{2}}$. Here ${{C}_{0}}$ is a constant independent of $\Gamma $.
Therefore, we have obtained (\ref{lb-est2}) and completed the proof of Theorem~\ref{theomIV}.
\medskip

\section{Concluding Remarks}\label{CR}
\medskip

The {\em virial theorem} in physics provides a relationship between the time-average of the total kinetic energy and that of the potential energy. For quantum multi-particle systems governed by the linear Schr\"{o}dinger equation, this often results in an elegant ratio. However, for the Schr\"{o}dinger equation in optics with non power-law type nonlinearities such as those square-root and saturable types, no virial results were available previously, to the best of our knowledge. Our study has yielded results concerning the virial relation and also the energy estimate of the {\em ground state}. Still, not too much is known about the higher energy states.
\medskip

\noindent {\bf Appendices: Some Technical Propositions and Lemmas}
\medskip

\noindent {\bf Appendix~I}
\medskip

\noindent {\bf Proposition A.1.}~~Let $f\left( s \right)=\frac{\left( 1-\frac{1}{1+s} \right)s}{s-\ln \left( 1+s \right)}$ for $s>0$. Then ${f}'\left( s \right)<0$ for $s>0$, $\underset{s\to 0+}{\mathop{\lim }}\,f\left( s \right)=2$ and $\underset{s\to \infty }{\mathop{\lim }}\,f\left( s \right)=1$.
\begin{proof}
By direct calculation,
\begin{eqnarray*}
{f}'\left( s \right) &=& {{\left( 1+s \right)}^{-2}}{{\left[ s-\ln \left( 1+s \right) \right]}^{-2}}\left[ 2{{s}^{2}}-\left( {{s}^{2}}+2s \right)\ln \left( 1+s \right) \right] \\
 &=& {{\left( 1+s \right)}^{-2}}{{\left[ s-\ln \left( 1+s \right) \right]}^{-2}}\rho \left( s \right) \,, \\
\end{eqnarray*}
for $s>0$, where $\rho \left( s \right)=2{{s}^{2}}-\left( {{s}^{2}}+2s \right)\ln \left( 1+s \right)$. We claim that $\rho \left( s \right)<0$ for $s>0$. By direct calculation, ${\rho }'\left( s \right)=s\eta \left( s \right)$ and $\eta \left( s \right)=4-2\frac{1+s}{s}\ln \left( 1+s \right)-\frac{2+s}{1+s}$ for $s>0$. Moreover, ${\eta }'\left( s \right)=-\frac{2}{{{s}^{2}}}\omega \left( s \right)$, $\omega \left( s \right)=s-\ln \left( 1+s \right)-\frac{{{s}^{2}}}{2{{\left( 1+s \right)}^{2}}}$ and ${\omega }'\left( s \right)=\frac{{{s}^{2}}}{{{\left( 1+s \right)}^{2}}}+\frac{{{s}^{2}}}{{{\left( 1+s \right)}^{3}}}>0$  for $s>0$. Then ${\eta }'\left( s \right)<0$ for $s>0$, which gives $\eta \left( s \right)<\eta \left( 0 \right)=0$ for $s>0$. Hence ${\rho }'\left( s \right)=s\eta \left( s \right)<0$ for $s>0$, which implies $\rho \left( s \right)<\rho \left( 0 \right)=0$ and ${f}'\left( s \right)={{\left( 1+s \right)}^{-2}}{{\left[ s-\ln \left( 1+s \right) \right]}^{-2}}\rho \left( s \right)<0$ for $s>0$ as well. On the other hand, we easily get $\underset{s\to 0+}{\mathop{\lim }}\,f\left( s \right)=2$ and $\underset{s\to \infty }{\mathop{\lim }}\,f\left( s \right)=1$ by direct calculation. Therefore, we have completed the proof of Proposition~A.1.
\end{proof}
\medskip

\noindent{\bf Appendix~II. Proof of Theorem~B}
\medskip

In order to be able to apply the argument of the proof of Theorem~A (cf.~\cite {LC-jmp-2014}), we remark that ${{\left( \sqrt{1+{{u}^{2}}}-1 \right)}^{2}}\le \frac{1}{4}{{u}^{4}}$ for $u\in \mathbb{R}$ (due to the fact that $0\le \sqrt{1+{{u}^{2}}}-1\le \frac{1}{2}{{u}^{2}}$ for $u\in \mathbb{R}$), which is almost same as the crucial inequality ${{u}^{2}}-\ln \left( 1+{{u}^{2}} \right)\le \frac{1}{2}{{u}^{4}}$ for $u\in \mathbb{R}$ in~\cite{LC-jmp-2014}. Then as for Proposition 3.1 and 3.2 in~\cite{LC-jmp-2014}, we have the analogs as folllows.

\noindent{\bf Proposition~A.2.}~~{\it
\begin{enumerate}
\item[(i)]~~Suppose $\Gamma \in \left( 0,{{T}_{3}} \right)$, i.e., $0<\Gamma <{{T}_{3}}$. Then the value ${{\tilde{e}}_{\Gamma }}=0$ can not be attained by a minimizer.
\item[(ii)]~Suppose $\Gamma \le 0$. Then ${{\tilde{e}}_{\Gamma }}=0$ and the value can not be attained by a minimizer.
\end{enumerate}
}

\noindent Hence, we can complete the proof of Theorem~B~(i) by combining Propositions~A.2 and~A.3.

For the proof of Theorem~B~(ii), we use the same ideas as in~\cite{LC-jmp-2014} and consider the following problem:
$$
{{\tilde{e}}_{\Gamma ,\varepsilon }}=\underset{\begin{matrix}
   u\in H_{0}^{1}\left( {{B}_{\frac{1}{\varepsilon }}} \right)  \\
   {{P}_{\varepsilon }}\left[ u \right]=1  \\
\end{matrix}}{\mathop{\inf }}\,{{\tilde{E}}_{\Gamma ,\varepsilon }}\left[ u \right]\,,
$$
where ${{B}_{\frac{1}{\varepsilon }}}$is the ball with radius $\frac{1}{\varepsilon }$ centered at the origin in $\mathbb{R}^2$,
\[{{\tilde{E}}_{\Gamma ,\varepsilon }}\left[ u \right]=\frac{1}{2}\int_{{{B}_{\frac{1}{\varepsilon }}}}{{{\left| \nabla u \right|}^{2}}-\Gamma {{\left( \sqrt{1+{{u}^{2}}}-1 \right)}^{2}}}\,, \]
${{P}_{\varepsilon }}\left[ u \right]=\int_{{{B}_{\frac{1}{\varepsilon }}}}{{{u}^{2}}}$ for $\varepsilon >0$ and $u\in H_{0}^{1}\left( {{B}_{\frac{1}{\varepsilon }}} \right)$. Note that ${{\left( \sqrt{1+{{u}^{2}}}-1 \right)}^{2}}={{u}^{2}}+2\left( 1-\sqrt{1+{{u}^{2}}} \right)\le {{u}^{2}}$ and $\Gamma >{{T}_{3}}>0$ so the lower bound estimate ${{\tilde{E}}_{\Gamma ,\varepsilon }}\left[ u \right]\ge -\frac{1}{2}\Gamma $ holds true. Hence, we may apply a symmetric-decreasing rearrangement (see Chapter~3 in~\cite{LL}) and the truncation argument (see the proof of Lemma 3.3 in~\cite{LC-jmp-2014}) to prove the following lemmas.
\medskip

\noindent {\bf Lemma~A.4.}~~{\it Assume $\Gamma >{{T}_{3}}>0$.
\begin{enumerate}
\item[(i)]~~For $\varepsilon >0$,  ${{\tilde{e}}_{\Gamma ,\varepsilon }}$ can be achieved by a minimizer ${{u}_{\varepsilon }}={{u}_{\varepsilon }}\left( r \right)\geq 0$ which is a function radially symmetric and monotone decreasing with $r$.
\item[(ii)]~For $\varepsilon >0$ sufficiently small, ${{\tilde{e}}_{\Gamma ,\varepsilon }}\le -{{c}_{0}}$ , where ${{c}_{0}}$ is a positive constant independent of~$\varepsilon $.
\end{enumerate}
}
\medskip

\noindent {\bf Lemma~A.5}~~{\it Under the same hypothesis of Lemma~A.4, minimizer ${{u}_{\varepsilon }}$ satisfies $${{\left\| {{u}_{\varepsilon }} \right\|}_{{{H}^{1}}\left( {{B}_{\tfrac{1}{\varepsilon }}} \right)}}\le {{K}_{0}}\,,$$ for $\varepsilon >0$ sufficiently small, where $K_0$ is a positive constant independent of $\varepsilon $. }
\medskip

The minimizer ${{u}_{\varepsilon }}$ satisfies the following equation:
	
\BE\label{eq-mizr1}
\Delta {{u}_{\varepsilon }}+\Gamma \left( 1-\frac{1}{\sqrt{1+u_{\varepsilon }^{2}}} \right){{u}_{\varepsilon }}={{\lambda }_{\varepsilon }}{{u}_{\varepsilon }}\text{   in  }{{B}_{\tfrac{1}{\varepsilon }}}\,,
\EE
with the zero Dirichlet boundary condition ${{u}_{\varepsilon }}=0\text{  on  }\partial {{B}_{\tfrac{1}{\varepsilon }}}$, where ${{\lambda }_{\varepsilon }}$ is the associated Lagrange multiplier. Multiply equation (\ref{eq-mizr1}) by ${{u}_{\varepsilon }}$ and integrate over ${{B}_{\tfrac{1}{\varepsilon }}}$. Then using integration by parts and $\int_{{{B}_{\tfrac{1}{\varepsilon }}}}{u_{\varepsilon }^{2}}=1$, we get ${{\lambda }_{\varepsilon }}={{\lambda }_{\varepsilon }}\int_{{{B}_{\tfrac{1}{\varepsilon }}}}{u_{\varepsilon }^{2}}=-{{\int_{{{B}_{\tfrac{1}{\varepsilon }}}}{\left| \nabla {{u}_{\varepsilon }} \right|}}^{2}}+\Gamma \int_{{{B}_{\frac{1}{\varepsilon }}}}{\left( 1-\frac{1}{\sqrt{1+u_{\varepsilon }^{2}}} \right)}u_{\varepsilon }^{2}$. Hence, by Lemma~A.5,
\BE\label{est-ev}
\left| {{\lambda }_{\varepsilon }} \right|\le {{K}_{1}},
\EE
where $K_1$ is a positive constant independent of $\varepsilon $. Here we have used the fact that $0\le \left( 1-\frac{1}{\sqrt{1+u_{\varepsilon }^{2}}} \right)u_{\varepsilon }^{2}\le u_{\varepsilon }^{2}$ and $\|u_\varepsilon\|_2=1$. Since ${{u}_{\varepsilon }}={{u}_{\varepsilon }}\left( r \right)$ is radially symmetric, equation (\ref{eq-mizr1}) and the zero Dirichlet boundary condition can be reduced to a boundary value problem of an ordinary differential equation as follows:
	
\BE\label{eq-mizr2}
\left\{ \begin{matrix}
   {{{{u}''}}_{\varepsilon }}+\frac{1}{r}{{{{u}'}}_{\varepsilon }}+\Gamma \left( 1-\frac{1}{\sqrt{1+u_{\varepsilon }^{2}}} \right){{u}_{\varepsilon }}={{\lambda }_{\varepsilon }}{{u}_{\varepsilon }}\text{  for  }0<r<\tfrac{1}{\varepsilon },  \\
   {{{{u}'}}_{\varepsilon }}\left( 0 \right)=0,\text{                    }{{u}_{\varepsilon }}\left( \tfrac{1}{\varepsilon } \right)=0.  \\
\end{matrix} \right.
\EE
Then by the uniqueness of ordinary differential equations, we can show the following.
\medskip

\noindent {\bf Lemma~A.6.}~~{\it The minimizer ${{u}_{\varepsilon }}={{u}_{\varepsilon }}\left( r \right)$ is positive for $0<r<\frac{1}{\varepsilon}$.}
\medskip

We may extend ${{u}_{\varepsilon }}$ to the entire plane ${{\mathbb{R}}^{2}}$ by setting ${{u}_{\varepsilon }}\left( r \right)=0\text{  for  }r>\tfrac{1}{\varepsilon }$. Note that each ${{u}_{\varepsilon }}$ is radially symmetric. Then Lemma~A.5 gives ${{\left\| {{u}_{\varepsilon }} \right\|}_{H_{r}^{1}\left( {{\mathbb{R}}^{2}} \right)}}\le {{K}_{0}}$, yielding
\BE\label{conv1}
{{u}_{\varepsilon }}\to U\text{  weakly in   }H_{r}^{1}\left( {{\mathbb{R}}^{2}} \right)
\EE
as $\varepsilon $ goes to zero (up to a subsequence).
Furthermore, by the compact embedding of $H^1_r$ radial functions in the space of $L^4_r$ functions (cf.~\cite{L-fa82}), we have
\BE\label{conv2}
{{u}_{\varepsilon }}\to U\text{  in   }L_{r}^{4}\left( {{\mathbb{R}}^{2}} \right)
\EE
as $\varepsilon $ goes to zero (up to a subsequence).
Now, we want to prove that $U$ is nontrivial. Due to ${{\left( \sqrt{1+{{u}^{2}}}-1 \right)}^{2}}\le \frac{1}{4}{{u}^{4}}$ for $u\in \mathbb{R}$, we get ${{\left( \sqrt{1+u_{\varepsilon }^{2}}-1 \right)}^{2}}\le \frac{1}{4}u_{\varepsilon }^{4}$. Hence, by (\ref{conv1}), (\ref{conv2}) and Lemma~A.4~(ii),
\begin{eqnarray*}
-{{c}_{0}} &\ge& \underset{\varepsilon \to 0+}{\liminf}\,{{\tilde{e}}_{\Gamma ,\varepsilon }} \\
&=&\underset{\varepsilon \to 0+}{\liminf}\,\frac{1}{2}{{\int_{{{\mathbb{R}}^{2}}}{\left| \nabla {{u}_{\varepsilon }} \right|^2}}}-\Gamma {{\left( \sqrt{1+u_{\varepsilon }^{2}}-1 \right)}^{2}}, \\
&\ge& \underset{\varepsilon \to 0+}{\liminf}\,\frac{1}{2}{{\int_{{{\mathbb{R}}^{2}}}{\left| \nabla {{u}_{\varepsilon }} \right|^2}}}-\tfrac{1}{4}\Gamma u_{\varepsilon }^{4} \\
&\ge& \frac{1}{2}{{\int_{{{\mathbb{R}}^{2}}}{\left| \nabla U \right|^2}}}-\tfrac{1}{4}\Gamma {{U}^{4}}\,,
\end{eqnarray*}
showing that $U$ is nontrivial. Here Fatou's Lemma has been used. Otherwise, $U\equiv 0$ and we would have $0>-c_0\geq 0$, a contradiction.

Now we claim that the limit function $U$ satisfies
	
\BE\label{eq-U1}
\Delta U+\Gamma \left( 1-\frac{1}{\sqrt{1+{{U}^{2}}}} \right)U={{\lambda }_{0}}U\text{   in   }{{\mathbb{R}}^{2}}\,,
\EE
$U=U\left( r \right)\ge 0$ is radially symmetric, and $\underset{r\to \infty }{\mathop{\lim }}\,U\left( r \right)=0$, where ${{\lambda }_{0}}$ is the limit of ${{\lambda }_{\varepsilon }}$'s (up to a subsequence) since (\ref{est-ev}) implies
\BE\label{conv-ev}
{{\lambda }_{\varepsilon }}\to {{\lambda }_{0}}\text{   as   }\varepsilon \to \text{0+  (up to a subsequence)}.
\EE
Let $\phi \in C_{0}^{\infty }\left( {{\mathbb{R}}^{2}} \right)$ be any test function. Since ${{u}_{\varepsilon }}$ satisfies (\ref{eq-mizr1}), we have
\BE\label{weq-minzr1}
\int_{{{\mathbb{R}}^{2}}}{\nabla {{u}_{\varepsilon }}\cdot \nabla \phi -\Gamma \int_{{{\mathbb{R}}^{2}}}{\left( 1-\frac{1}{\sqrt{1+u_{\varepsilon }^{2}}} \right){{u}_{\varepsilon }}\phi =-{{\lambda }_{\varepsilon }}\int_{{{\mathbb{R}}^{2}}}{{{u}_{\varepsilon }}\phi }}}\,.
\EE
Hence, (\ref{conv1}) and (\ref{conv2}) give
$$\int_{{{\mathbb{R}}^{2}}}{\nabla {{u}_{\varepsilon }}\cdot \nabla \phi \to \int_{{{\mathbb{R}}^{2}}}{\nabla U\cdot \nabla \phi }}\,, $$
$$\int_{{{\mathbb{R}}^{2}}}{{{u}_{\varepsilon }}\phi \to \int_{{{\mathbb{R}}^{2}}}{U\phi }}\,, $$
and
\begin{eqnarray*}
& & \int_{{{\mathbb{R}}^{2}}}{\frac{{{u}_{\varepsilon }}}{\sqrt{1+u_{\varepsilon }^{2}}}} \phi =\int_{{{\mathbb{R}}^{2}}}{\left( {{u}_{\varepsilon }}-U \right)\frac{\phi }{\sqrt{1+u_{\varepsilon }^{2}}}}+\int_{{{\mathbb{R}}^{2}}}{\frac{1}{\sqrt{1+u_{\varepsilon }^{2}}}}U\phi  \\
& & =\int_{{{\mathbb{R}}^{2}}}{\left( {{u}_{\varepsilon }}-U \right)\frac{\phi }{\sqrt{1+u_{\varepsilon }^{2}}}}+\int_{{{\mathbb{R}}^{2}}}{\frac{1}{\sqrt{1+{{U}^{2}}}}}U\phi +\int_{{{\mathbb{R}}^{2}}}{\left( \frac{1}{\sqrt{1+u_{\varepsilon }^{2}}}-\frac{1}{\sqrt{1+{{U}^{2}}}} \right)}U\phi  \\
& & =\int_{{{\mathbb{R}}^{2}}}{\left( {{u}_{\varepsilon }}-U \right)\frac{\phi }{\sqrt{1+u_{\varepsilon }^{2}}}}+\int_{{{\mathbb{R}}^{2}}}{\frac{1}{\sqrt{1+{{U}^{2}}}}}U\phi \\
& & \hspace{0.5cm} +\int_{{{\mathbb{R}}^{2}}}{\left( U-{{u}_{\varepsilon }} \right)\frac{U+{{u}_{\varepsilon }}}{\sqrt{1+u_{\varepsilon }^{2}}\sqrt{1+{{U}^{2}}}\left( \sqrt{1+{{U}^{2}}}+\sqrt{1+u_{\varepsilon }^{2}} \right)}} U\phi \\
& & \to \int_{{{\mathbb{R}}^{2}}}{\frac{1}{\sqrt{1+{{U}^{2}}}}}U\phi\,.
\end{eqnarray*}
Note that ${{\left\| \frac{\phi }{\sqrt{1+u_{\varepsilon }^{2}}} \right\|}_{\tfrac{4}{3}}}\le {{\left\| \phi  \right\|}_{\tfrac{4}{3}}}$ and ${{\left\| \frac{U+{{u}_{\varepsilon }}}{\sqrt{1+u_{\varepsilon }^{2}}\sqrt{1+{{U}^{2}}}\left( \sqrt{1+{{U}^{2}}}+\sqrt{1+u_{\varepsilon }^{2}} \right)}U\phi  \right\|}_{\frac{4}{3}}}\le {{\left\| U\phi  \right\|}_{\frac{4}{3}}}\le {{\left\| U \right\|}_{4}}{{\left\| \phi  \right\|}_{2}}$ as
$$
0\le \frac{U+{{u}_{\varepsilon }}}{\sqrt{1+u_{\varepsilon }^{2}}\sqrt{1+{{U}^{2}}}\left( \sqrt{1+{{U}^{2}}}+\sqrt{1+u_{\varepsilon }^{2}} \right)}\le 1\,.
$$
Thus, (\ref{weq-minzr1}) and (\ref{conv-ev}) imply
$$
\int_{{{\mathbb{R}}^{2}}}{\nabla U\cdot \nabla \phi }-\Gamma \int_{{{\mathbb{R}}^{2}}}{\left( 1-\frac{1}{\sqrt{1+{{U}^{2}}}} \right)U}\phi =-{{\lambda }_{0}}\int_{{{\mathbb{R}}^{2}}}{U\phi }\,,
$$
from which $U$ satisfies (\ref{eq-U1}) and $\underset{\left| x \right|\to \infty }{\mathop{\lim }}\,U\left( x \right)=0$. Moreover, $U=U\left( r \right)\ge 0$ is radially symmetric, since each ${{u}_{\varepsilon }}$ is positive and radially symmetric. Therefore, the equation for $U$ can be written as
$$
\left\{ \begin{matrix}
   {U}''+\frac{1}{r}{U}'+\Gamma \left( 1-\frac{1}{\sqrt{1+{{U}^{2}}}} \right)U={{\lambda }_{0}}U\text{   for   }r>0,  \\
   {U}'\left( 0 \right)=0,\text{              }U\left( \infty  \right)\text{=0                 }.  \\
\end{matrix} \right.
$$
By the uniqueness of ordinary differential equations, we have obtained the following.
\medskip

\noindent{\bf Lemma~A.7}~~{\it $U\left( r \right)>0\text{   for   }r\ge 0$. \hspace{4cm} $\Box$ }
\medskip

Due to the fact that $\underset{r\to \infty }{\mathop{\lim }}\,U\left( r \right)=0$, there exists $R_1>0$ such that $0<U\left( r \right)\le 1\text{   for   }r\ge {{R}_{1}}$. By equation (\ref{eq-U1}), $\Delta U=\left( {{\lambda }_{0}}-\Gamma \left( 1-\frac{1}{\sqrt{1+{{U}^{2}}}} \right) \right)U\in {{L}^{4}}\left( {{B}_{{{R}_{1}}}} \right)$ , since ${{\lambda }_{0}}-\Gamma \left( 1-\frac{1}{\sqrt{1+{{U}^{2}}}} \right)\in {{L}^{\infty }}$. Hence, by the standard regularity theorem of the Poisson equation, $U\in {{W}^{2,4}}\left( {{B}_{{{R}_{1}}}} \right)$ and then by the Sobolev embedding ${{W}^{2,4}}\left( {{B}_{{{R}_{1}}}} \right)\subset {{L}^{\infty }}\left( {{B}_{{{R}_{1}}}} \right)$, we have obtained the following.
\medskip

\noindent {\bf Lemma~A.8}~~{\it $U\left( r \right)\le {{K}_{2}}\text{   for  }r\ge 0$, where $K_2$ is a positive constant.\hspace{4cm} $\Box$ }
\medskip

Now we prove that ${{\lambda }_{0}}$ is positive by contradiction. Suppose ${{\lambda }_{0}}\le 0$. Then equation (\ref{eq-U1}) and Lemma~A.7 imply
\[\Delta U={{\lambda }_{0}}U-\Gamma \left( 1-\frac{1}{\sqrt{1+{{U}^{2}}}} \right)U\le 0\text{   in  }{{\mathbb{R}}^{2}}\,.\]
Therefore, by Lemma~A.8 and the Liouville Theorem, $U$ must be a constant function, i.e., $U\equiv 0$, which is impossible. Therefore, we conclude the following.
\medskip

\noindent {\bf Lemma~A.9.}~~{\it The limit ${{\lambda }_{0}}$ of ${{\lambda }_{\varepsilon}}$ as $\varepsilon\downarrow 0$ satisfies ${{\lambda }_{0}}>0$. \hspace{6cm} $\Box$ }
\medskip

By (\ref{conv2}), $u_\varepsilon$ converges to $U$ almost everywhere as $\varepsilon\downarrow 0$ (up to a subsequence). Moreover, since each $u_\varepsilon$ is monotone decreasing with $r$, we have the following.
\medskip

\noindent {\bf Lemma~A.10.}~~{\it $U=U\left( r \right)$ is monotone decreasing with $r$. \hspace{6cm} $\Box$}
\medskip

To complete the proof of Theorem~B, we now only need to prove that ${{\left\| U \right\|}_{2}}=1$ if ${{u}_{\varepsilon }}\to U$ strongly in ${{L}^{2}}\left( {{\mathbb{R}}^{2}} \right)$ as $\varepsilon \to 0$ (up to a subsequence). Note that each $u_{\varepsilon }$ satisfies ${{\left\| {{u}_{\varepsilon }} \right\|}_{2}}=1$. Fix $\sigma >0$ and consider the set ${{N}_{\sigma ,\varepsilon }}=\left\{ r>0:{{u}_{\varepsilon }}\left( r \right)>\sigma  \right\}$. Then by Lemmas~A.4~(i) and~A.6, ${{N}_{\sigma ,\varepsilon }}=\left( 0,{{R}_{\sigma ,\varepsilon }} \right)$ and
\[\pi {{\sigma }^{2}}R_{\sigma ,\varepsilon }^{2}\le 2\pi \int_{0}^{{{R}_{\sigma ,\varepsilon }}}{ru_{\varepsilon }^{2}\left( r \right)dr}\le 2\pi \int_{0}^{\infty }{ru_{\varepsilon }^{2}\left( r \right)dr}=\left\| {{u}_{\varepsilon }} \right\|_{2}^{2}=1,\]
implying
\BE\label{upb-Rsigma}
{{R}_{\sigma ,\varepsilon }}\le \frac{1}{\sigma \sqrt{\pi }}.
\EE
Note that the upper bound of ${{R}_{\sigma ,\varepsilon }}$ is $\frac{1}{\sigma \sqrt{\pi }}$ independent of $\varepsilon $. Hence, (\ref{conv1}) gives
$$
{{u}_{\varepsilon }}\to U\text{  weakly in   }H_{r}^{1}\left( {{B}_{\tfrac{1}{\sigma \sqrt{\pi }}}} \right)\,,
$$ and the standard Sobolev compact embedding ${{H}^{1}}\left( {{B}_{\tfrac{1}{\sigma \sqrt{\pi }}}} \right)\hookrightarrow {{L}^{2}}\left( {{B}_{\tfrac{1}{\sigma \sqrt{\pi }}}} \right)$ implies
\BE\label{conv1-1}
{{u}_{\varepsilon }}\to U\quad\text{  strongly in   }\:{{L}^{2}}\left( {{B}_{\tfrac{1}{\sigma \sqrt{\pi }}}} \right)
\EE
as $\varepsilon $ goes to zero (up to a subsequence). Moreover, (\ref{upb-Rsigma}) implies $r\notin {{N}_{\sigma ,\varepsilon }}$ and ${{u}_{\varepsilon }}\left( r \right)\le \sigma $ for $r>\frac{1}{\sigma \sqrt{\pi }}$ . Let $\sigma >0$ be sufficiently small such that $\frac{{{\lambda }_{0}}}{2}>\Gamma {{\sigma }^{2}}$. Then by (\ref{eq-mizr2}), (\ref{conv-ev}) and Lemma~A.9,
\begin{eqnarray*}
{{{{u}''}}_{\varepsilon }}+\frac{1}{r}{{{{u}'}}_{\varepsilon }}
&=&-\Gamma \left( 1-\frac{1}{\sqrt{1+u_{\varepsilon }^{2}}} \right){{u}_{\varepsilon }}+{{\lambda }_{\varepsilon }}{{u}_{\varepsilon }}  \\
&\ge& \left( -\Gamma {{\sigma }^{2}}+{{\lambda }_{\varepsilon }} \right){{u}_{\varepsilon }}  \\
&\ge& \frac{{{\lambda }_{0}}}{2}{{u}_{\varepsilon }}\quad\quad\quad\text{   for   }\frac{1}{\sigma \sqrt{\pi }}<r<\frac{1}{\varepsilon }\,,  \\
\end{eqnarray*}
i.e.,
\BE\label{eq-mizr2-1}
{{{u}''}_{\varepsilon }}+\frac{1}{r}{{{u}'}_{\varepsilon }}\ge \frac{{{\lambda }_{0}}}{2}{{u}_{\varepsilon }}\quad\text{   for   }\quad\frac{1}{\sigma \sqrt{\pi }}<r<\frac{1}{\varepsilon }\,.
\EE
Here we have used the fact that $1-\frac{1}{\sqrt{1+u_{\varepsilon }^{2}}}=\frac{u_{\varepsilon }^{2}}{\sqrt{1+u_{\varepsilon }^{2}}\left( \sqrt{1+u_{\varepsilon }^{2}}+1 \right)}\le u_{\varepsilon }^{2}\le {{\sigma }^{2}}$ for $r>\frac{1}{\sigma \sqrt{\pi }}$. Note that ${{\lambda }_{0}}$ is a positive constant independent of $\varepsilon$. Consequently, equation (\ref{eq-mizr2-1}) gives
\BE\label{conv1-2}
{{u}_{\varepsilon }}\left( r \right)\le {{e}^{-\alpha r}}\quad\text{   for   }\quad r>\frac{1}{\sigma \sqrt{\pi }},
\EE
where $\alpha >0$ is a constant independent of $\varepsilon $. Therefore, by (\ref{conv1-1}), (\ref{conv1-2}) and the Dominated Convergence Theorem, we get
${{u}_{\varepsilon }}\to U$ strongly in ${{L}^{2}}\left( {{\mathbb{R}}^{2}} \right)$, so the proof of Theorem~B is complete.

\noindent \textbf{Acknowledgement}. This work was initiated when T.C. Lin was visiting Texas A and M University at Qatar (TAMUQ) in 2013 and finalized during his visit in 2016. He would like to thank TAMUQ for kind hospitality during his visits. The research of T.C. Lin was also partially supported by the National Center for Theoretical Sciences (NCTS) and NSC grant 103-2115-M-002-005 of Taiwan. The work of other coauthors has been partially supported by the Qatar
National Research Fund project NPRP 8-028-1-001.

\end{document}